\documentclass{amsart}
\usepackage{graphicx}    
\usepackage{amssymb,mathrsfs,amscd}
\usepackage{Mytheorems}
\author{Alexander Sotirov}
\address{Stanford University} 
\author{Shih-Hsien Yu}
\address{City University of Hong Kong}
\thanks{The research of the second author is supported by CERG grant CityU 103304}

\dedicatory{Dedicated to  Professor Tai-Ping Liu on the occasion of his 
60th birthday}

\begin{document}     
\title{On the Boltzmann Diffusion of Two Gases}

\maketitle
\begin{abstract}

We study the Boltzmann equation for a mixture of two gases in one space dimension
 with initial condition of one
 gas near vacuum and the other near a Maxwellian equilibrium state.
 A qualitative-quantitative mathematical analysis is developed to study this mass diffusion 
problem based on the Green's function of the Boltzmann equation for the single species 
hard sphere collision model in \cite{liuyug1}. The cross species resonance of the
mass diffusion and the diffusion-sound wave is investigated. An 
exponentially sharp global  solution is obtained.

\end{abstract} 
 
\begin{section}{Introduction}
The Boltzmann equation for a mixture of two gases in one space dimension takes the form of the
 following system:
$$
 \begin{cases} \displaystyle  
\partial_t{ F_A}  + \xi^1 \partial_x { F_A} =  Q^{AA}(F_A, F_A) + Q^{AB}(F_A,F_B),
 \\ \\
\displaystyle  
\partial_t { F_B}  + \xi^1 \partial_x{F_B}=  Q^{BB}(F_B, F_B) + Q^{BA}(F_B,F_A).
\end{cases}
$$
The right hand side consists of the usual collision terms which for $X,Y \in\{A,B\}$ are given by, \cite{bird}:
$$ Q^{XY}(F_X, F_Y)= \sigma_{XY}^2 \int_{S^+}\int_{\mathbb{R}^3} (F_X'{F_Y}'_* - F_X{F_Y}_* )
 |(\xi-\xi_*) \cdot n| d \xi_* dn$$
where ${F_Y'}_*= F_Y(x,t, \xi'_*)$, ${F_Y}_*=F_Y(x,t,\xi_*)$, and $F_X'=F_X(x,t,\xi')$, with $S$
 the unit sphere and $S^+$ are $n$ such that $(\xi-\xi_*)\cdot n >0$. The constant factor
 $\sigma_{XY}$ denotes the sum of the radii of molecule $X$ and molecule $Y$. The
 post-collision velocities $\xi'$ and $\xi_*'$ are given by the formulas:
$$ \begin{cases} \displaystyle  
\xi'=\xi + \frac {2m_Y}{m_X+m_Y}((\xi_*-\xi) \cdot n) n, \\ \\
\displaystyle  
\xi_*'=\xi_* - \frac {2m_X}{m_X+m_Y}((\xi_*-\xi) \cdot n) n,
\end{cases}
$$
where $m_A$ and $m_B$ are the molecular masses of the two gases. These formulas are obtained from the
 conservation of momentum and energy at collision and the assumption that forces only act normally to 
the surface.

The studies on gas mixture in terms of 
 Boltzmann equations  were initiated
 in \cite{takata2} to study the 
condensation-vaporization problem for  mixture  of 
vapors of different species. Many interesting 
physical problems such as ``ghost effect'' and 
``Knudsen layer'' for gas mixture have been investigated,
\cite{Takata,takata3}.

This research is also interested in phenomena related to vapor-vapor mixtures. Its focus 
is on the basic phenomenon of mass diffusion of a finite amount of gas $A$ 
dissolved into the surrounding gas $B$. From this consideration, we write 
$$\begin{cases}
F_A= {\mathsf f}_A {\mathsf M}_A^{\frac{1}{2}} \text{ for gas with  finite total mass amount} , \\
F_B= {\mathsf f}_B {\mathsf M}_B^{\frac{1}{2}} + {\mathsf M}_B \text{ for gas in the background}.
\end{cases} 
$$
Here, ${\mathsf M}_A$ and ${\mathsf M}_B$ are the Maxwellian states

$$ \begin{cases} \displaystyle  
  {\mathsf M}_B(\xi) \equiv   \frac {\rho_B}{\sqrt { (2\pi  RT)^3 } } \exp 
\left({ \frac {-|\xi|^2}{2RT}} \right), \\
\displaystyle  {\mathsf M}_A(\xi) \equiv {\mathsf M}_B
\left(\frac{ \sqrt{m_A} }{ \sqrt{m_B}    } \xi \right)=
\frac {1}{\sqrt { (2\pi  RT)^3 } } \exp \left( - \frac{ m_A }{ m_B   } { \frac {|\xi|^2}{2RT}} \right).
\end{cases}
$$

Here $T$ is the temperature of the mixture, $\rho_B$ the density of gas $B$, $R$ the gas constant, and we
 have assumed that the average velocity is zero.
We also assume that 
$$ RT =1.$$
The two Maxwellian states satisfy 
\begin{equation}
\label{ABBA}
 \begin{cases}
Q^{AA}({\mathsf M}_A,{\mathsf M}_A)=0, \; \;  Q^{BB}({\mathsf M}_B,{\mathsf M}_B)=0,\\
Q^{AB}({\mathsf M}_A,{\mathsf M}_B)=0, \; \;Q^{BA}({\mathsf M}_B,{\mathsf M}_A)=0.
\end{cases} 
\end{equation}

A mixture of rarefied gases is important for many high technology industrial applications.  A classical
 approach to the study of gas mixtures is through a Navier-Stokes model. However, to use the Navier-Stokes
 model one needs to obtain the diffusion coefficient through formal calculations of a transport equation
 which are carried out easily only for either small or large ratio of the molecular masses \cite{landau}. This
 makes the foundation of the Navier-Stokes approach uncertain, especially since the Navier-Stokes equation
 is valid for fluid near thermo-equilibrium state. Often in applications the state of the gases is not
 necessarily close to  such a thermo-equilibrium. The validity of the Navier-Stokes approach is also
 questionable if there is some physical boundary encountered. 
Such issues on validity for Navier-Stokes in the application to 
rarefied gases  had been extensively studied by an engineering school in 
Kyoto University initiated by Professor Sone, \cite{sone}. 
Thus, it is desirable to study the problem
 entirely through a fundamental physical model based on a kinetic equation.

The system for $({\mathsf f}_A, {\mathsf f}_B)$ is
$$
\begin{cases}
\partial_t {\mathsf f}_A + \xi^1 \partial_x {\mathsf f}_A = L_{AB} {\mathsf f}_A + \Gamma^{AA}({\mathsf f}_A,{\mathsf f}_A)
+
\Gamma^{AB}({\mathsf f}_A, {\mathsf f}_B), \\
\partial_t {\mathsf f}_B + \xi^1 \partial_x {\mathsf f}_B = L  {\mathsf f}_B
+ \Gamma^{BB}({\mathsf f}_B, {\mathsf f}_B) +
  L_{BA} {\mathsf f}_A+  \Gamma^{BA}({\mathsf f}_B,
{\mathsf f}_A),
\end{cases} 
$$
where 
\begin{equation}
\label{n1.5}
\begin{cases} \displaystyle  L_{AB} {\mathsf f_A}_ \equiv \frac{ 1 }{  \sqrt{{\mathsf M}_A}  } Q^{AB}({\mathsf f}_A \sqrt{{\mathsf M}}_A, {\mathsf M}_B)  \\ 
\displaystyle  
 L_{BA} {\mathsf f}_A \equiv \frac{ 1 }{  \sqrt{{\mathsf M}_B}  } 
Q^{BA}( {\mathsf M}_B, {\mathsf f}_A \sqrt{{\mathsf M}_A} ),
\\
\displaystyle  
L {\mathsf f}_B \equiv   \frac{ 1 }{  \sqrt{{\mathsf M}_B}  }
\left(  Q^{BB}({\mathsf f}_B \sqrt{{\mathsf M}}_B, {\mathsf M}_B) +
 Q^{BB}( {\mathsf M}_B, {\mathsf f}_B \sqrt{{\mathsf M}}_B) \right), \\ \displaystyle  
\Gamma^{XY}({\mathsf f}_X,{\mathsf f}_Y) \equiv 
\frac{ 1 }{ \sqrt{{\mathsf M}_X}    } Q^{XY}({\mathsf f}_X \sqrt{{\mathsf M}_X}, {\mathsf f}_Y 
\sqrt{{\mathsf M}_Y}) \text{ for } X,Y \in \{ A,B\}.
\end{cases}
\end{equation}
Here, $\Gamma^{XY}$ are all quadratic nonlinear terms. 

When  both the total mass of gas $A$ and the perturbations in gas $B$ are sufficiently small, in general 
the basic time asymptotic behavior 
of the nonlinear solution is governed by a linear equation (up to a large time scale). It is reasonable
 that scent diffusion in a gas is a fast time scale
phenomenon independently of whether the density of scent is strong or weak. Therefore in great generality
 the following decoupled linear Boltzmann system 

\begin{align}
&\partial_t {\mathsf g} + \xi^1 \partial_x {\mathsf g} = L_{AB} {\mathsf g},  \label{lseq1}\\
&\partial_t {\mathsf h} + \xi^1 \partial_x {\mathsf h} =L {\mathsf h} +    L_{BA} {\mathsf g}, \label{lseq2}
\end{align}
will give a  refined  quantitative description of mass diffusion. It can be written in a matrix form as
\begin{equation}
\label{upper}
\partial_t \begin{pmatrix}
{\mathsf g} \\ {\mathsf h} 
\end{pmatrix} + \begin{pmatrix}
\xi^1 \partial_x - L_{AB} & 0 \\
- L_{BA} & \xi^1 \partial_x - L   
\end{pmatrix}  \begin{pmatrix}
{\mathsf g} \\ {\mathsf h} 
\end{pmatrix} = 0.
\end{equation} 
We will call \eqref{lseq1} the linear Boltzmann diffusion equation, and we will call the system \eqref{upper}
 the linear Boltzmann system.
There are several features which provide a qualitative difference between the problem at hand and the single
 gas flow or flow in a mixture where both gases are near a nonzero equilibrium state:

First, a purely diffusive behavior is dominant in the gas near vacuum and so we have qualitatively recovered
 the behavior already expected from a classical fluid treatment. Acoustic waves are not present in the first gas, such waves however appear in the second gas.

Second, as discussed in the paper \cite{Takata} , the technique that would be natural in dealing with the
 problem of two gases, both near nonzero background states, would require to work in a functional space of 
pairs of functions $({\mathsf f}_A,{\mathsf f}_B)$ with an inner product that depends on the densities
 $\rho_A,\rho_B$:
$$\langle({\mathsf f}_A,{\mathsf f}_B),({\mathsf h}_A,{\mathsf h}_B)\rangle=\rho_A\langle {\mathsf f}_A,
{\mathsf h}_A\rangle +\rho_B\langle {\mathsf f}_B, {\mathsf h}_B\rangle.$$
This inner product will yield a self-adjoint operator on the space of pairs of functions on which a
 spectral analysis can be performed and a theory similar to that developed in \cite{liuyug1}
 can be constructed. Notice however that if one of the gases is near vacuum (so $\rho_A=0$) the inner
 product becomes degenerate and it becomes impossible to work with the operator on pairs of functions. 
This is the problem that is addressed by this work.

Third, the characteristics for the  
macroscopic equations in the system  \eqref{lseq1}, \eqref{lseq2}
 coincide along the direction of mass diffusion.  A further physical property of
the collision operator, namely a microscopic cancellation representing the conservation of the total mass is employed to resolve this resonance. 

To see the generic feature of mass diffusion the following form of the initial data 
$({\mathsf g},{\mathsf h})|_{t=0} =({\mathsf g}_{in}, {\mathsf h}_{in})$  will be sufficient:

\begin{equation}
 \label{indata}
\begin{cases}
\displaystyle  {\mathsf g}_{in}(x, \xi)= {\mathsf h}_{in}(x,\xi) \equiv 0 \text{ for } |x| \ge 1, \\
\displaystyle  \sup_{ \begin{subarray}c |x| \le 1  \\ 
 \xi \in {\mathbb R}^3 
  \end{subarray}
} (1+|\xi|)^3\max\left(|{\mathsf h}_{in}(x,\xi)|, |{\mathsf h}_{in}(x,\xi)|\right)
\le 1.
\end{cases} 
\end{equation}

\begin{theorem}[Main Theorem] \label{maintheolin}
Suppose the initial data are as in \eqref{indata}. Then the solution to \eqref{lseq1} , \eqref{lseq2} 
satisfies the estimates \begin{eqnarray*}
\|{\mathsf g}(x,t) \|_{L^2_{\xi}} & = &  O(1)\left[ 
\frac{ e^{ - \frac{ (x- \lambda_2 t)^2 }{C(1+t)    } }  }{  \sqrt{(1+t) } } 
 + e^{-(t+|x|)/C} \right], \\
\|{\mathsf h}(x,t) \|_{L^2_{\xi}} & = & \ O(1)\left[ \sum_{i=1}^3
\frac{ e^{ - \frac{ (x- \lambda_i t)^2 }{C(1+t)    } }  }{  \sqrt{(1+t) } } 
 + e^{-(t+|x|)/C} \right] + O(1)
\begin{cases}
0 \text{ for } |x| \ge 2 |\lambda_1|t, \\
e^{-t^{\frac{1}{2}/C}} \text{ for } |x| \le
2 |\lambda_1 t|,
\end{cases} 
\end{eqnarray*}
\end{theorem}
where $C>0$,  $$ \begin{cases}
\lambda_1=- \sqrt{\frac{ 5 }{ 3   }  }, \; \; 
\lambda_2=0, \; \; \lambda_3=\sqrt{\frac{ 5  }{ 3   }  }, \\ \displaystyle  
\|{\mathsf f}(x,t) \|_{L^2_{\xi}} = 
\left(\int_{{\mathbb R}^3}{\mathsf f}(x,t,\xi)^2 d\xi\right)^{\frac{1}{2}}.
\end{cases}
$$

Within a finite Mach number region $|x| \le 2 |\lambda_1|(1+t)$, there is 
a tail decaying  exponentially in time: $O(1) e^{- t^{\frac{1}{2}}/C}$. This is
a purely kinetic phenomenon  (such a tail should not appear in a  
linearized equation modeled properly after the Navier-Stokes equation). 
However, since this tail decays exponentially in time, it is difficult to 
recognize the difference in the solutions of kinetic equations and classic
continuum fluid mechanics when the time scale is large.

The basic principle used to develop the theory is
 the {\bf ``separation of scales''}, which is a concept created in 
\cite{liuyug1}. In \cite{liuyug1}, a Long Wave-Short Wave decomposition and a Particle-Wave decomposition were introduced. These two decompositions are 
designed for solutions with different physical characteristics. For example, 
when the solution is in a region corresponding to the ``small time scale'' or the ``high
Mach number region'' (large space scale), the  behavior is particle-like. In this case the Particle-Wave decomposition is used to determine the dominant behavior.  When the solution is in
a region corresponding to large time scale or in the finite Mach number region, 
the behavior of the solution is rather close to continuum flow. In this case Long Wave-Short Wave decomposition is used to analyze the solution. These two decompositions can be matched to each other with the difference decaying exponentially fast in time in the norm  $\| \cdot \|_{L^\infty_x(L^2_\xi)}$  Note there is no exponentially decaying structure in the $x$ variable. This
is an instance of the uncertainty principle from harmonic analysis.
The Long Wave-Short Wave decomposition allows one to work in Fourier space and there use the microscopic cancellation to resolve the resonance problem due to the waves along the characteristic $x=0$.

In Section II we will describe the physical problem clearly and 
review some of the basic theory of the Boltzmann equation 
for a single gas relying on \cite{pinsky,liuyug1}.
In Sections III and IV  we establish the spectral analysis for  the Boltzmann
Diffusion equation and continue to establish its Green's function.

In Section V,  a microscopic cancellation  and a {\bf ``scale separation''}
 are discovered for the gas mixture problem. The main theorem is proved.

We would like to note that the above result allows an extension to a solution of the fully nonlinear equation
 provided the initial data is sufficiently close to equilibrium. The first author has carried out such an
 extension as a part of his dissertation in the case when the molecular mass satisfy the relation $m_A>m_B$.
 While the proof involves lengthy computations and a number of genuinely nontrivial difficulties had to be
 overcome, the main techniques are similar to those used for the nonlinear problem in \cite{liuyug1}.
\\ \\
{\bf Acknowledgment.}
The authors would like to thank Professor Tai-Ping Liu. His penetrating vision on both 
scientific and mathematical problems has been influencing the second author from the very beginning of his scientific career. Both authors  have greatly benefited from his teaching, his sharing of scientific visions and ideas. They  wish to express their gratitude to Professor Tai-Ping Liu on the occasion of his sixtieth birthday, the first author being a current graduate student of his, and the second one of his former graduated students,  a friend, and a long time collaborator.

\end{section}

\begin{section}{Preliminaries}
\begin{subsection}{A Physics Problem from Rarefied Gases} \hfill \\
The primary content of this problem is to study the  diffusion of  Gas A with finite amount of mass
in Gas B which occupies the entire space of  a one-dimensional pipe in terms of a system of Boltzmann equations. 
Since this is a linear coupled system, the linear mass diffusion of Gas A will 
have some resonance from the background particle diffusion of Gas B. The
linear resonance may also have a significant effect to the acoustic waves 
propagating in Gas B. It would be interesting to determine which one of  Figure 1
and  Figure 2 below gives a correct picture of the mass  diffusion for rarefied gases (notice the difference in spatial decay in the two figures).
 \\
\\ 

 \setlength{\unitlength}{1cm}
\begin{picture}(12,4.5)
\put(-1.5,0.5){\includegraphics[height=4cm]{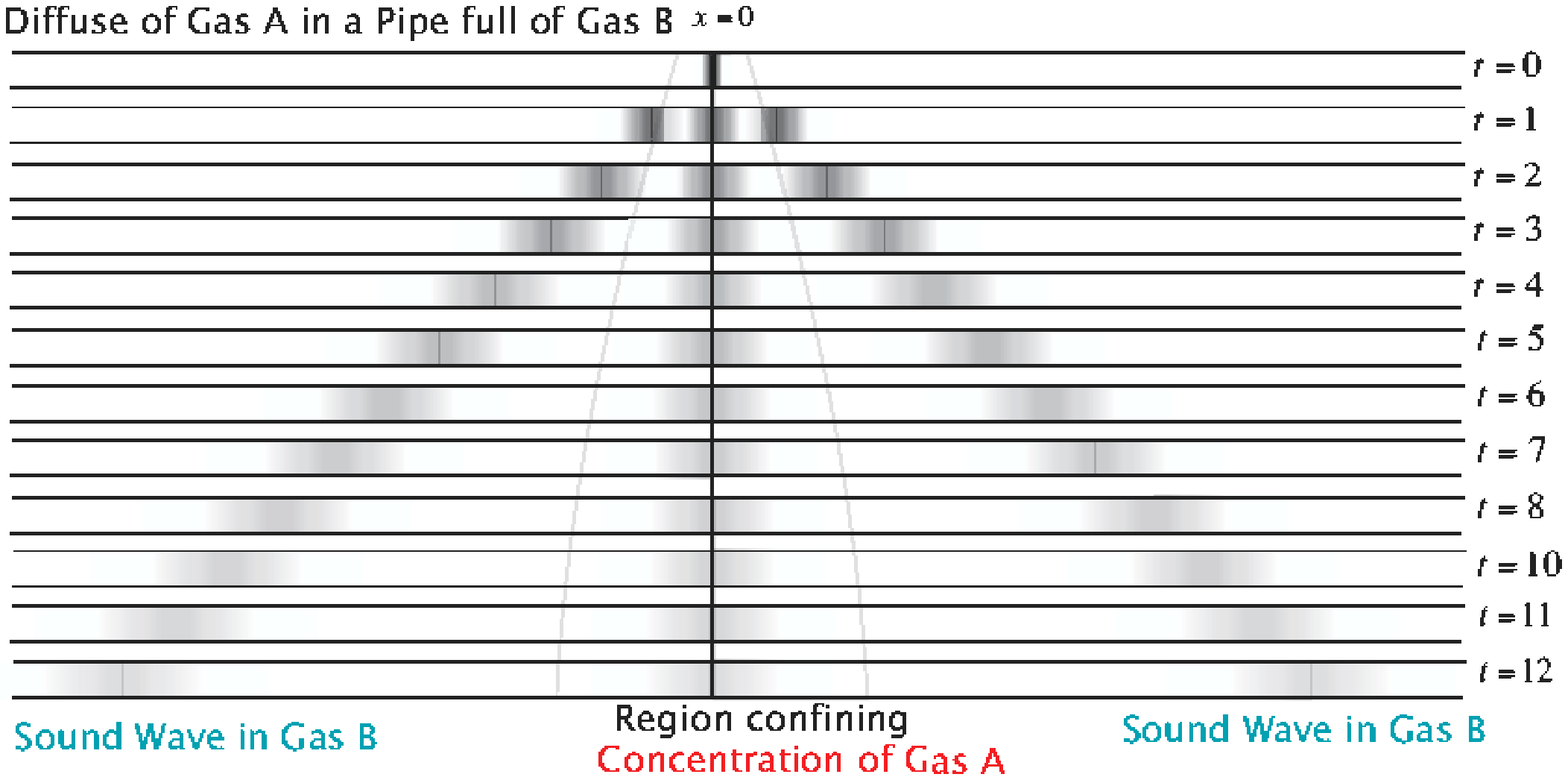}}
\put(7,0.5){\includegraphics[height=4cm]{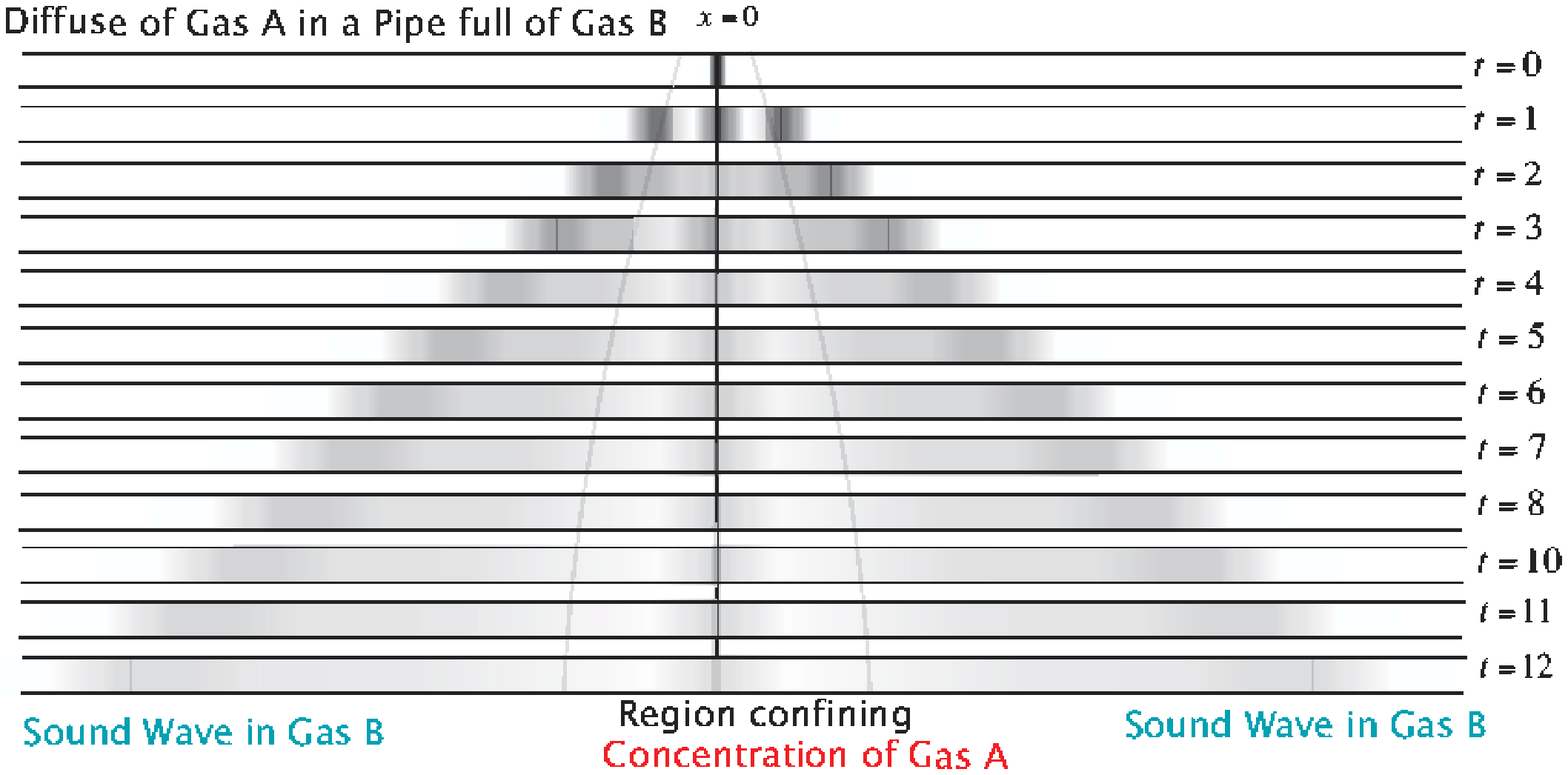}}
\put(-0.3,0.2){Figure 1.}
\put(8.2,0.2){Figure 2.}
\end{picture}

Here, Figure 1 is a picture based on common experience about gases close to
an equilibrium state. However, our knowledge for a gas near thermo-equilibrium
  is not necessary valid for a rarefied gas.  Figure 2 is a picture adopted from
the  consideration of resonance and coupling of two different gases. One might conjecture that there is a tail algebraic in space after the sound wave.

In general in physics, one may need to set up a physical experiment
to  judge which figure is close to the reality of a rarefied gas. 
Unfortunately,  to setup an experiment for this problem could not be easy in
practice. 
Since the question is about the coupling of mass diffusion to sound wave
propagation in a rarefied gas, it would require a rarefied gas environment with a  
large length scale.  

Another alternative in physics available nowadays is a numerical simulation from first principles which
 may give some information about the truth. 
 However, the requirement for a larger  space size  and the number of particles\footnote{ 
The particle density of a rarefied gas at $10^{-6}$ atm, $273 {}^\circ K$ is in the order 
$10^{16}/\ell$.
} makes a reliable numerical simulation very difficult.

Thus, a quantitative and qualitative mathematical analysis would be a good alternative for this
problem. Such a mathematical analysis would benefit other physics problems 
in the field of rarefied gases. 
\end{subsection} 

\begin{subsection}{Mathematical Theories on Boltzmann equation} \hfill \\

For the convenience of the reader,  in this section we include the necessary
 theorems 
established for a linear Boltzmann equation:
\begin{equation}
\label{lboe}
\partial_t {\mathsf b} + \xi^1 \partial_x {\mathsf b} =L {\mathsf b},
\end{equation} 
where $L$ is the linear collision operator given in \eqref{n1.5}.
The linearized collision operator $L$ is symmetric on the student Hilbert space 
$L^2_\xi$ given by
$$ \begin{cases}
\displaystyle  ({\mathsf h}, {\mathsf g}) \equiv
\int_{{\mathbb R}^3} {\mathsf h}(\xi) {\mathsf g}(\xi)
d \xi \text{ for } {\mathsf h}, {\mathsf g} \in
L^2_\xi, \\
\displaystyle  \|{\mathsf h} \|_{L^2_\xi}= \sqrt{
({\mathsf h}, {\mathsf h}) }.
\end{cases}
$$
For a given function ${\mathsf g}(x,t,\xi)$, we have the following
notation:
$$ \begin{array}{||c|c|c||} \hline \text{Expression} & \text{Definition} & \text{Comment} \\ \hline
\| {\mathsf g} \|_{L^2_\xi} & \displaystyle  \left(
\int_{{\mathbb R}^3} {\mathsf g}(x,t,\xi)^2 d \xi \right)^{\frac{1}{2}} & \text{ non-negative real-valued
 function in $(x,t)$} \\
\hline \|{\mathsf g} \|_{L^2_x(L^2_\xi)} & \displaystyle  
\left(\int_{\mathbb R} \int_{{\mathbb R}^3} 
{\mathsf g}(x,t,\xi)^2 d \xi dx \right)^{\frac{1}{2}} & \text{ non-negative real-valued function in $t$}
\\ \hline 
\|{\mathsf g} \|_{H^k_x(L^2_\xi)} & \displaystyle  
\left( 
\int_{\mathbb R} \int_{{\mathbb R}^3} \sum_{l=0}^k |\partial_x^l{\mathsf g}(x,t,\xi)|^2 d 
\xi dx \right)^{\frac{1}{2}} & \text{ non-negative real-valued function in $t$} \\
 \hline 
\|{\mathsf g} \|_{L^\infty_x(L^2_\xi)} & \displaystyle  
\sup_{  x \in  {\mathbb R}}
\|{\mathsf g}\|_{L^2_\xi} & \text{ non-negative real-valued function in $t$}
\\ \hline 
\|| {\mathsf g} \|| & \displaystyle  
\sup_{ \begin{subarray}c  x \in {\mathbb R} \\
\xi \in {\mathbb R}^3 \end{subarray}} |{\mathsf g}(x,t,\xi)|
(1+|\xi|)^3 & \text{ non-negative real-valued function in $t$}
\\
\hline 
\end{array} 
$$

For the linearized collision operator $L$ linearized around 
${\mathsf M} \equiv {\mathsf M}_B$, the null space of ${L}$ is a five-dimensional vector space with 
orthogonal basis
$\chi_i$, $i=0,1,2,3,4$, \cite{boltzmann}: 

$$
\begin{array}{l}
\displaystyle    ker({L}) \equiv span \{ \chi_0, \chi_1, \chi_2, \chi_3, \chi_4 \},  \\ 
\displaystyle   \begin{cases} 
\chi_0 \equiv {\mathsf M}^{1/2} \\ 
\chi_i \equiv \xi^i {\mathsf M}^{1/2} \text{ for } i=1,2,3, \\ 
\chi_4 \equiv \frac{1}{\sqrt{6}} (|\xi|^2 -3) {\mathsf M}^{1/2}. 
\end{cases}  
\end{array}
$$
  
Since we are interested in   planar wave propagations,  we  restrict the 
 function space to: 
\begin{equation}
 \label{***} 
  \displaystyle   
L^2_\xi  \equiv \left\{ {\mathsf g} \in L^\infty({\mathbb R}^3)|  
(\xi^2 {\mathsf M}^{1/2},{\mathsf g})=0, \; 
(\xi^3 {\mathsf M}^{1/2}, {\mathsf g})=0, \; \| {\mathsf g}  \|_{L^2_\xi} <
\infty \right\},  \; \;
  \xi \equiv (\xi^1,\xi^2,\xi^3). 
\end{equation}

\noindent{\bf Macro-Micro Decomposition} \\
This space is invariant under the full Boltzmann operator 
$Q^{BB}$  and 
the linearized operator $L$. We decompose the Hilbert space 
$L_\xi^2 = \ker({L}) \oplus \ker({L})^\perp:$  
for any ${\mathsf g} \in L_\xi^2$,  
$$ \begin{cases} 
 {\mathsf g} \equiv {\mathsf P}_0 {\mathsf g} + {\mathsf P}_1 {\mathsf g} 
( \equiv {\mathsf g}_0 + {\mathsf g}_1), \\ 
{\mathsf P}_0 {\mathsf g} \equiv  
 ( \chi_0, {\mathsf g} ) \chi_0 +  
( \chi_1, {\mathsf g} ) \chi_1 +  
( \chi_4 , {\mathsf g} ) \chi_4, \\ 
{\mathsf P}_1 {\mathsf g} \equiv {\mathsf g} - {\mathsf P}_0 {\mathsf g}. 
\end{cases}  
$$ 
The characteristic information of the Euler equations is  
connected to the operator ${\mathsf P}_0 \xi^1$ on ${\mathsf P}_0 L_{\xi}^2$, \cite{liuyu}: 
\begin{equation} 
\label{MM2} 
\begin{split} 
& dim( {\mathsf P}_0 L_\xi^2)=3, \\ 
 &{\mathsf P}_0 \xi^1 {\mathsf E}_i = \lambda_i {\mathsf E}_i  
\text{ for } i=1,2,3,  \\ 
&\left\{ \lambda_1 = - \sqrt{5/3}, \;  
\lambda_2=0, \; \lambda_3 = \sqrt{5/3} \right\}, \\ 
&\begin{cases}{\mathsf E}_1 \equiv \left( 
\sqrt{\frac{3  }{ 2   }}\chi_0 - \sqrt{ \frac{5  }{ 2   }}\chi_1 + \chi_4 \right), \\  
{\mathsf E}_2 \equiv \left(- \sqrt{\frac{2  }{3    } }\chi_0 + \chi_4 \right) , \\ 
{\mathsf E}_3 \equiv \left(  
\sqrt{\frac{ 3 }{ 2   } }\chi_0 + \sqrt{ \frac{5  }{2    } }\chi_1 + \chi_4 
\right), \\ 
({\mathsf E}_i,{\mathsf E}_j)=\delta_j^i \text{ (Kronecker's delta function).} 
\end{cases}  
\end{split}  
\end{equation}  
 $$  {\boldsymbol c} \equiv |\lambda_1|:
\text{ Speed of sound at rest.}  $$

For the hard sphere model we consider, the linearized collision operator $L$ is of the following form, \cite{hilbert}: 
 
\begin{equation} 
\label{In3} 
 \begin{cases} \begin{array}{rl} 
L {\mathsf g} (\xi) &= -\nu(\xi) {\mathsf g}(\xi) + {\mathsf K} {\mathsf g} (\xi),  
\\ 
\displaystyle  {\mathsf K} {\mathsf g} (\xi) &\equiv  
\int_{{\mathbb R}^3} K(\xi,\xi_*) {\mathsf g}(\xi_*) \; d \xi_*, \\ 
K(\xi,\xi_*)  &\equiv \displaystyle   
\frac{ 2 \sigma_{BB}^2   }{ \sqrt{2 \pi}|\xi -\xi_*|   } \exp\left(  - \frac{(|\xi|^2 - |\xi_*|^2)^2  }{8 | \xi - \xi_*|^2    }  
 - \frac{| \xi-\xi_*|^2  }{8    }  \right) \\ & \hfill \displaystyle   
- \frac{ |\xi -\xi_*| }{ 2   } \exp \left( 
- \frac{  (|\xi|^2 + |\xi_*|^2 )}{  4  }  
\right),  \\ 
\displaystyle  \nu(\xi) &  \displaystyle  \equiv 
\frac{\sigma^2_{BB}}{ \sqrt{2 \pi}    } \left( 2 e^{-\frac{|\xi|^2  }{  2  } } + 
2 \left( |\xi|+ \frac{ 1 }{ |\xi|   }  \right) \int_0^{|\xi|} e^{-\frac{ u^2  }{ 2   } } \; du  
\right) \ge \nu_0(1+|\xi|), \\
\nu(\xi) &  \displaystyle \sim 1+ | \xi|. 
\end{array}  
\end{cases}  
\end{equation}

  \begin{lemma}
  There exists positive constant $\nu_0$ such that, for any ${\mathsf h} \in L^2_\xi$,
$$
    \begin{cases}
\displaystyle  ({\mathsf P}_0 {\mathsf h}, {\mathsf P}_0 {\mathsf h})
      = \sum_{j=1}^3 ({\mathsf E}_j, {\mathsf h})^2, \\
      \displaystyle  
({\mathsf P}_0 {\mathsf h},\xi^1 {\mathsf P}_0 {\mathsf h})
    =\sum_{j=1}^3 \lambda_j ({\mathsf E}_j,{\mathsf h})^2,
    \\
    ({\mathsf P}_1 {\mathsf h}, L {\mathsf P}_1 {\mathsf h}
    ) \le - \nu_0 ({\mathsf P}_1 {\mathsf h},
    (1+|\xi|){\mathsf P}_1 {\mathsf h}).
    \end{cases}
    $$
  
    \label{L3.3}
\end{lemma}

  This lemma follows from direct computations and Carleman's
  theory on the negative definiteness of the operator $L$ on 
  $Range({\mathsf P}_1)$, \cite{carleman}.
\\ \\

\noindent{\bf Estimates on Green's function.} \cite{liuyug1} \\ 
 The Green's function ${\mathbb G}(x,t,\xi,\xi_*)$ for the initial value problem \eqref{lboe} constructed
  in \cite{liuyug1} plays a fundamental role for
 studying the Boltzmann  diffusion equation. \\
It has also been shown in \cite{liuyug1} that ${\mathbb G}(x,t,\xi,\xi_*)$ also
satisfies the backwards equation:
$$
\begin{cases} 
(- \partial_\tau - \xi_* \partial_y - L) {\mathbb G}
  (x-y,t-\tau,\xi,\xi_*)=0, \\
  {\mathbb G}(x-y,0,\xi,\xi_*) = \delta^1(y-x) \delta^3(
  \xi_* - \xi);
\end{cases}
$$
and the function 
$${\mathsf b}(x,t,\xi)\equiv 
\int_{\mathbb R} \int_{{\mathbb R}^3} {\mathbb G}(x-y,t,\xi,\xi_*)
{\mathsf b}(y,0,\xi_*) d \xi_* dy$$
solves the initial value problem \eqref{lboe}. 

Without spelling out the parameters
$\xi$ and $\xi_*$ in the Green's function ${\mathbb G}(x,t,\xi,\xi_*)$, one can treat ${\mathbb G}(x,t)$ as an
  $L_\xi^2$-operator-valued function in $(x,t)$:\\
  For given $(x,t) \in {\mathbb R}\times {\mathbb R}^+$,
$$
  \begin{array}{rl}
\displaystyle    {\mathbb G}(x,t):& {\mathsf h}  \in L^2_\xi \longmapsto
  {\mathbb G}(x,t) {\mathsf h} \in L^2_\xi, \\
\displaystyle  {\mathbb G}(x,t) {\mathsf h}(\xi) &
\displaystyle  \equiv
\int_{{\mathbb R}^3} {\mathbb G}(x,t,\xi,\xi_*) {\mathsf h}(
\xi_*) d \xi_*.
\end{array} 
  $$

\noindent{\bf Analytic Spectral Decomposition}

Take Fourier transformation of \eqref{lboe} in $x$ variable, then the
solution ${\mathsf b}(x,t)$ can be expressed 
$$
\hat {\mathsf b}(k,t) = e^{(-i k \xi^1 +L)t} \hat {\mathsf b}(k,0) 
\text{ for } k \in {\mathbb R}.
$$

\begin{lemma}[Ellis-Pinsky]
\hfill \\ Assume $L^2_\xi$ satisfies  the restriction \eqref{***}.\\
 \label{pinsky-liu} 
\noindent (i) For any $\kappa_0>0$ if $|k|>\kappa_0>0$, $k \in \mathbb{R}$ then there exists $ \epsilon(\kappa_0) > 0$ such that if $\sigma$ is
 an eigenvalue of $-ik\xi^1+L$ then $Re (\sigma) < - \epsilon$.  \\
\noindent (ii) There exist $\kappa_0$ and $\delta$ such that for $|k|<\kappa_0$ the spectrum with $|\sigma(k)|<\delta$ 
consists of three points $\sigma_1(k), \sigma_2(k), \sigma_3(k)$ which converge to the origin as $k$ tends zero and they , together with the corresponding eigenvectors are 
{\bf analytic } functions of $k$ for $|k|<\kappa_0$.  \\
\noindent (iii) More precisely, the expansions of the eigenvectors 
${\mathsf e}_j(k)$ and eigenvalues $\sigma_j(k)$  are as follows:

\begin{equation} 
\label{B1.5} 
 \begin{cases} \begin{array}{rl} 
\sigma_1(k) &=-i \lambda_1 \; k + A_1 |k|^2 + O(1) | k|^3,  \\ 
\sigma_2(k) &= -i \lambda_2 \; k + A_2 |k|^2 +O(1)|k|^3, \\ 
\sigma_3(k) &=-i \lambda_3 \; k + A_3 |k|^2 + O(1) | k|^3, \\ 
\displaystyle   
A_j &=  \left( {\mathsf P}_1 \xi^1 {\mathsf E}_j, {L}^{-1} {\mathsf P}_1 
\xi^1 {\mathsf E}_j \right)<0, 
\\ 
A_3&=A_1, \end{array}  
\end{cases} 
\end{equation}  

\begin{equation} 
\label{normal_} 
 \begin{cases} 
(-i \xi^1 k +{L}) {\mathsf e}_j(k) = \sigma_j(k) {\mathsf e}_j(k), \\  
\left({{\mathsf e}}_j(k),{\mathsf e}_k\right(k))=\delta_{jk}, \\ 
{\mathsf e}_i(k) = {\mathsf E}_i + k {\mathsf e}_i'(0) + O(1) |k|^2, 
\end{cases}   
\end{equation}  
and, 
$$ 
\begin{array}{l} \displaystyle   
{\mathsf e}_i'(0) = \sum_{j=1}^3 \varepsilon_i^j {\mathsf E}_j +  
{\mathsf e}_i^\perp, \; \;  {\mathsf e}_i^\perp  
\equiv {\mathsf P}_1 {\mathsf e}_i'(0), 
\\ 
\begin{cases} \displaystyle   
\varepsilon_k^j = -i \frac{ ({\mathsf E}_j, {\mathsf P}_0  
\xi^1 {L}^{-1} {\mathsf P}_1 \xi^1 {\mathsf E}_k)  }{ (\lambda_j - \lambda_k) }\text{ for } j \ne k, \\ 
\varepsilon_k^k=0 , \\ 
{{\mathsf e}}_k^\perp = i {L}^{-1} {\mathsf P}_1 \xi^1 {\mathsf E}_k. 
\end{cases}  
\end{array}  
$$ 
\end{lemma}
Note that the analyticity property of $\sigma_j(k)$ for $k$ around $0$ 
follows from theorems in \cite{kato}. 
\begin{lemma}[Liu-Yu]
There exists $\kappa_0>0$ such that for all
$|k| \le \kappa_0$ the operator  $e^{(-ik \xi^1 +L)t}$ on $L^2_\xi$ can be
decomposed as
\begin{equation}
\label{spec_d1}
 e^{(-ik \xi^1 +L)t} = \sum_{j=1}^3 e^{\sigma_j(k)t} {\mathsf e}_j(k)
\otimes \left \langle {\mathsf e}_j(k) \right| + e^{(-ik \xi^1 +L)t}
\Pi_k^\perp
\end{equation}
so that 
\begin{equation}
\label{spec_d2}
\left\| e^{(-ik \xi^1 +L)t} \Pi_k^\perp \right\|_{L^2_\xi} \le C e^{-t/C} \text{ for some } C>0.
\end{equation}

\end{lemma} 
\noindent Here, the operator $ {\mathsf e}_j(k)
\otimes \left \langle {\mathsf e}_j(k) \right|$ is given as follows. \\
For any ${\mathsf j}$ and ${\mathsf k}$ in $L^2_\xi$, the operator
${\mathsf j} \otimes \left \langle  {\mathsf k}  \right| $ on
$L^2_\xi$ is defined as follows:
$$
 {\mathsf j} \otimes \left \langle  {\mathsf k}  \right| {\mathsf g}
\equiv ({\mathsf k}, {\mathsf g}) {\mathsf j}.
$$

 Now taking inverse Fourier transform in the generalized sense and writing as a convolution 
we get for the solution:
$$
{\mathsf g} (x,t)  =\mathbb{G}^t {\mathsf g}_{in} 
\equiv \frac {1}{\sqrt{2\pi}} \int_{\mathbb{R}}e^{ikx+(-i\xi^1 k + L)t}\hat
{\mathsf g}_{in}(k) dk  
=\int_{\mathbb{R}}\mathbb{G}(x-y,t) {\mathsf g}_{in}(y) dy,
$$
where
$$ \mathbb{G}(x,t) \equiv \frac {1}{\sqrt{2\pi}} \int_{\mathbb{R}}e^{ikx+(-i\xi^1 k + L)t} dk$$
and we accept the last identity as the definition of Green's function $\mathbb{G}(x,t)$.
$\mathbb{G}(x,t)$ is a generalized operator valued function of $x$ with values operators on
$L^2_\xi=L^2(\mathbb{R}^3(\xi))$ since we obtained it by taking inverse Fourier transform of the
regular operator valued function $e^{(-i\xi^1k+L)t}$ which as a function of $k$ does not necessarily
 have a convergent Fourier integral. By pointwise estimates we mean estimates of its operator norm as 
a function of the points $x$ and $t$: but since the function is generalized of course one needs to 
identify and subtract its singularity for such estimates to make sense. However at this stage we will 
only look for estimates once Green's function is convoluted with a bounded function with compact 
support in the $x$ variable, so pointwise estimates make sense; i.e. we will look to estimate
 $\|\mathbb{G}g_{in}\|$ where $g_{in}$ is as in \eqref{indata}.

\noindent{\bf Long Wave-Short Wave Decomposition.} \cite{liuyug1} \\
One can decompose the semi-group 
${\mathbb G}(x,t) \equiv \displaystyle  \frac{ 1 }{ \sqrt{2\pi}    }  \int_{{\mathbb R}} e^{ix k + (-i k \xi^1 +L)t} d k$ into 
\begin{equation}
\label{l-s-de}
\begin{cases} \displaystyle  
 {\mathbb G}_L(x,t) \equiv   \frac{ 1 }{ \sqrt{2\pi}    }
 \int_{|k| \le \kappa_0 } e^{ix k + (-i k \xi^1 +L)t} d k,
 \\ \displaystyle  
 {\mathbb G}_S(x,t) \equiv    \frac{ 1 }{ \sqrt{2\pi}    } \int_{|k|\ge \kappa_0} e^{ix k + (-i k \xi^1 +L)t} d k.
\end{cases} 
\end{equation}

The Green's function is obtained in terms of the Fourier variable $k$ without detailed spectral properties of the operator $(-i k \xi^1 +L)$ described when $|k|\ge \kappa_0$. Thus, it not possible to obtain $\| {\mathbb G}_S(x,t)\|_{L^2_\xi}$ pointwise through Fourier analysis alone. This is a matter of uncertainty principle.

{\bf Particle-Wave Decomposition.} \cite{liuyug1} 
\\
Consider an initial value problem 
\begin{equation}
\label{p-w-d}
\begin{cases}
\partial_t {\mathsf b} + \xi^1 \partial_x {\mathsf b} = L {\mathsf b}, \\
{\mathsf b}(x,0)= {\mathsf g}_{in},
\end{cases} 
\end{equation}
where  ${\mathsf g}_{in}$ satisfies
\begin{equation}
\label{gin}
\begin{cases}
{\mathsf g}_{in}(x,\xi) \equiv 0 \text{ for } |x| \ge 1, \\
\|| {\mathsf g}_{in}
\|| \equiv \| {\mathsf g}_{in} \|_{L^\infty_x(L^\infty_{\xi,3})}=1
  < \infty.
\end{cases}
\end{equation}
\begin{lemma}[Liu-Yu] For each given $j \in {\mathbb N}$, 
there exist functions ${\mathbb P}_j(x,t)$ and ${\mathbb W}_j(x,t)$ so that
the solution ${\mathsf b}(x,t)$ of \eqref{p-w-d} with the given initial data in \eqref{gin} 
satisfies
$$ \begin{cases}
{\mathsf b}(x,t) = {\mathbb P}_j(x,t) + {\mathbb W}_j(x,t) \text{ for all }
x \in {\mathbb R}, t\ge 0, \\
{\mathbb P}_j(x,0) \equiv {\mathsf g}_{in}(x), \\
\|{\mathbb P}_j(x,t) \|_{L^2_\xi} \le C_j e^{-(|x|+t)/C_j} \text{ for all } x 
\in {\mathbb R}, t \ge 0, \\
\| (\partial_t + \xi^1 \partial_x - L){\mathbb P}_j \|_{H^j_x(L^2_\xi)} 
\le C_j e^{-t/C_j}, 
\\
\|{\mathbb W}_j\|_{H^j_x(L^2_\xi)} \le C_j \text{ for all } t \ge 0 \text{(uniformly bounded in time)},
\end{cases} 
$$
where $C_j>0$ is an universal constant.
\label{lemma-p-w-d}

\end{lemma} 
\begin{remark}
The existence of such a decomposition is due to both the introduction of an
essential kinetic equation (4.2) in \cite{liuyug1}  to resolve the singularity
in the initial data into source terms which are smooth in $\xi$ and to 
the discovery of Mixture Lemma in \cite{liuyug1}. 
\end{remark} 

\begin{theorem}[Liu-Yu]
\label{re-on-g}
For the $L^2_\xi$-operator-valued function ${\mathbb G}(x,t)$, Lemma \ref{lemma-p-w-d} is still valid, i.e. 
there exist  ${\mathbb P}_j(x,t)$ and ${\mathbb W}_j(x,t)$ such that
\begin{equation}
\label{re-on-1}
\begin{cases}
{\mathbb G}(x,t) = {\mathbb P}_j(x,t) + {\mathbb W}_j(x,t) \text{ for all } (x,t) \in {\mathbb R} \times {\mathbb R}^+, \\
{\mathbb P}_j(x,t) =  e^{-\nu(\xi_*)t} \delta(x-\xi^1 t) \delta^3(\xi-\xi_*)+ {\mathsf j}_1(x,t)+
{\mathsf j}_2(x,t)+ \cdots + {\mathsf j}_j(x,t), \\
\|{\mathsf j}_l(x,t) \|_{L^2_\xi} \le C_j e^{-(|x|+t)/C_j} \text{ for all } x \in {\mathbb R},\; \; l =2,3,\cdots,j,
\in {\mathbb R}, t \ge 0, \\
\| (\partial_t + \xi^1 \partial_x - L){\mathbb P}_j \|_{H^j_x(L^2_\xi)} 
\le C_j e^{-t/C_j}, 
\\
\|{\mathbb W}_j\|_{H^j_x(L^2_\xi)} \le C_j \text{ for all } t \ge 0 \text{(uniformly bounded in time)},
\end{cases} 
\end{equation} 
\end{theorem}

 We recall the following theorems in \cite{liuyug1} on the Green's function ${\mathbb G}$
 for the initial value problem.

\begin{theorem}[Liu-Yu] 
  \label{mainA} For the bounded initial data ${\mathsf g}_{in}$  
with compact support in $x$, \eqref{gin},  there exist 
positive constants $C_0$ and $C_1$ independent of ${\mathsf g}_{in}$ 
such that 
\begin{equation} 
\label{maint} 
\begin{split} 
&\text{for $|x |  \le 2  {\boldsymbol c}  t$,} \\ 
&\hskip0.5cm \begin{cases} \begin{array}{l} 
\displaystyle  \| {\mathbb G}^t {\mathsf g}_{in}(x)\|_{L^2_\xi} \displaystyle  =O(1)\; ||| {\mathsf g}_{in} ||| 
\left(  \sum_{i=1}^3 \frac{ e^{ - \frac{|x - \lambda_i t|^2  }{C_0 (t+1) }  }  } {\sqrt{ (t+1)    }}   
+ e^{-(|x|+t)/C_1} 
\right), \\ 
\displaystyle  \| {\mathbb G}^t {\mathsf P}_1 {\mathsf g}_{in}(x)\|_{L^2_\xi}, \;  
 \|  {\mathsf P}_1 {\mathbb G}^t {\mathsf g}_{in}(x)\|_{L^2_\xi} \displaystyle  \\  \displaystyle  
\hspace{2cm} =O(1)\; ||| {\mathsf g}_{in} ||| \left(  \sum_{i=1}^3
 \frac{ e^{ - \frac{|x - \lambda_i t|^2  }{C_0 (t+1) }  }  } {{ (t+1)    }}  + e^{-(|x|+t)/C_1} \right), \\ 
\displaystyle  \|  {\mathsf P}_1{\mathbb G}^t {\mathsf P}_1 {\mathsf g}_{in}(x)\|_{L^2_\xi}  
\displaystyle  =O(1) \; ||| {\mathsf g}_{in} ||| \left(  \sum_{i=1}^3 
\frac{ e^{ - \frac{|x - \lambda_i t|^2  }{C_0 (t+1) }  }  } {{ (t+1)^{3/2}    }}  + e^{-(|x|+t)/C_1} \right), \\ 
\displaystyle  \| {\mathbb G}^t {\mathsf E}_k \otimes \left \langle
 {\mathsf E}_k \right| {\mathsf g}_{in}(x)\|_{L^2_\xi} \\ \hspace{2cm}  \displaystyle  =O(1) \; ||| {\mathsf g}_{in} ||| 
\left(   \frac{ e^{ - \frac{|x - \lambda_k t|^2  }{C_0 (t+1) }  }  } {\sqrt{ (t+1)    }}  + 
 \sum \limits_{ \begin{subarray}c  1 \le i \le 3 \\ 
i \ne k \end{subarray}} \frac{ e^{ - \frac{|x - \lambda_i t|^2  }{C_0 (t+1) }  }  } {{ (t+1)    }}   
+ e^{-(|x|+t)/C_1} 
\right), \\ 
\end{array}  
\end{cases} \\ 
&\text{ for $|x| \ge 2  {\boldsymbol c}  t$} \\ 
&\hskip3.5cm  \|  {\mathbb G}^t {\mathsf g}_{in}(x) \|_{L_\xi^2} \le C\; ||| {\mathsf g}_{in} ||| 
\; e^{-(|x|+t)/C_1},
\end{split}  
\end{equation}
where 
  \begin{align}
\nonumber 
& 
  {\mathbb G}^t {\mathsf g}_{in}(x) \equiv
  \int_{{\mathbb R}} {\mathbb G}(x-y,t) {\mathsf g}_{in}(y)
  dy.
  \end{align}

\end{theorem}
\noindent
{\bf Remark.} This theorem is stated as Theorem 5.6 in \cite{liuyug1}. 
The first work to employ complex analysis to obtain such exponentially sharp 
estimates is in \cite{liu-zeng} for  the  Green's function of compressible
Navier-Stokes equation.   {\hfill
  $\boxtimes$}

\begin{theorem}[Liu-Yu] \label{mainc}
 For the bounded, compact-supported initial data ${\mathsf g}_{in}$, \eqref{gin}, there exist positive  
constants $C_1$ and   $C_0$ such that
for  all $x \in {\mathbb R}$ the main fluid part satisfies 
\begin{equation}
  \label{mainc1} 
  \| [ {\mathsf E}_j \otimes \left \langle {\mathsf E}_j 
\right | \; {\mathbb G}^t \; {\mathsf E}_k \otimes \left \langle {\mathsf E}_k 
\right| \; {\mathsf g}_{in}](x) \|_{L^\infty_{\xi,3}}    =O(1) ||| {\mathsf g}_{in} |||  
\left( \delta_j^k \frac{ e^{  - \frac{ (x - \lambda_k t)^2  }{ C_0 (t+1)  }      }  }{  \sqrt{(1+t)  }}  
+ \sum_{i=1}^3 \frac{ e^{  - \frac{ (x - \lambda_i t)^2  }{ C_0 (t+1)  }      }  }{  {(1+t)  }} + e^{-(|x|+t)/C_1} 
\right).
\end{equation}
and the non-fluid parts have higher rate of decay in time:
$$
  \| [{\mathbb G}^t {\mathsf P}_1 {\mathsf g}_{in}(x) \|_{L^\infty_{\xi,3}}  =O(1) ||| {\mathsf g}_{in} |||  
\left(  \sum_{i=1}^3 \frac{ e^{  - \frac{ (x - \lambda_i t)^2  }{ C_0 (t+1)  }      }  }{  {(1+t)  }}  +
 e^{-(|x|+t)/C_1} 
\right), 
$$
$$
  \|    {\mathsf P}_1 {\mathbb G}^t {\mathsf g}_{in}(x) \|_{L^\infty_{\xi,3}}  =O(1) ||| {\mathsf g}_{in} |||  
\left(  \sum_{i=1}^3 \frac{ e^{  - \frac{ (x - \lambda_i t)^2  }{ C_0 (t+1)  }      }  }{  {(1+t)  }}  +
 e^{-(|x|+t)/C_1} 
\right),  
$$
$$
 \|  {\mathsf P}_1  {\mathbb G}^t {\mathsf P}_1 {\mathsf g}_{in}(x) \|_{L^\infty_{\xi,3}}  =O(1) ||| {\mathsf g}_{in} |||  
\left(  \sum_{i=1}^3 \frac{ e^{  - \frac{ (x - \lambda_i t)^2  }{ C_0 (t+1)  }      }  }{  {(1+t)^{3/2}  }}  + e^{-(|x|+t)/C_1} 
\right). 
$$

\end{theorem}  
\end{subsection}

\end{section}

\begin{section}{Spectral Properties of the Linear Boltzmann Diffusive Equation }
\begin{subsection}{On the Cross Species Linear Collision Operator} \hfill \\
Similarly to the linear single species collision operator, the operator $L^{AB}$ is  symmetric in $L^2_\xi$, however now  $dim(ker(L^{AB}))=1$.  One has the following lemma: 
\begin{lemma} \hfill \\ \label{lemma31}
i) The linear cross species  collision operator $L^{AB}$ is symmetric in $L^2_\xi$. \\
ii)  $ ker(L^{AB}) = \{ c {\mathsf M}_A^{\frac{1}{2}}| c \in {\mathbb R}  \}$.
\end{lemma} 
\begin{proof}
From the definitions of both $L^{AB}$ in \eqref{n1.5} and  
$Q^{AB}$, one has that for any ${\mathsf u}, {\mathsf v} \in
L^2_\xi$,

$$ L_{AB} {\mathsf u}  = \frac{1}{{\mathsf M}_A^{1/2}} \sigma^2_{AB} \int_{S^+}\int_{\mathbb{R}^3} 
({\mathsf u}' {{\mathsf M}_A^{\prime1/2}} {{\mathsf M}_B}'_* - {\mathsf u}
 {\mathsf M}_A^{1/2}{{\mathsf M}_B}_* ) |(\xi-\xi_*) \cdot n| d \xi_* dn$$ 
where ${{\mathsf M}_B'}_*= {\mathsf M}_B( \xi'_*)$, 
${{\mathsf M}_B}_*={\mathsf M}_B(\xi_*)$, and ${\mathsf u}'={\mathsf u}(\xi')$. 

Next notice that ${\mathsf M}_A(|\xi|)={\mathsf M}_B( ( \frac {m_A}{m_B} )^{1/2} \xi)$ due to the energy consrevation equation satisfies 
$${\mathsf M}'_A{{{\mathsf M}'}_B}_*={\mathsf M}_A {\mathsf M}_{B*}$$
which allows us to compute (with $(,) $ denoting the usual inner product in $\mathbb{R}^3_\xi$) 
after the change $\xi \rightarrow \xi'$, $\xi_* \rightarrow \xi_*'$
\begin{align*} 
( L_{AB}{\mathsf u},{\mathsf v}  )  &= \sigma^2_{AB}\int \frac{{\mathsf v}}{{\mathsf M}_A^{1/2}}
 \int \left( {\mathsf u}'{\mathsf M}_A^{\prime 1/2} {{\mathsf M}_B}'_* - {\mathsf u} 
{\mathsf M}_A^{1/2}{{\mathsf M}_B}_* \right) (|(\xi-\xi_*)\cdot n|) dn d\xi_* d\xi \\
 &= \frac{1}{2} \sigma^2_{AB}\int ( \frac{{\mathsf v} }{{\mathsf M}_A^{1/2}} 
-\frac{{\mathsf v}'}{{{\mathsf M}_A}^{\prime 1/2}}) ( {\mathsf u}'{\mathsf M}_A^{\prime 1/2} {{\mathsf M}_B}'_* 
- {\mathsf u} {\mathsf M}_A^{1/2}{{\mathsf M}_B}_* ) (|(\xi-\xi_*)\cdot n|) dn d\xi_* d\xi \\
 &= -\frac{1}{2} \sigma^2_{AB}\int \frac{1}{{\mathsf M}_A{\mathsf M}_{B*}}
 ( {\mathsf v}'{{\mathsf M}_A}^{\prime 1/2}{\mathsf M}_{B*}'
 -{\mathsf v}{\mathsf M}_A^{1/2}{\mathsf M}_{B*}) 
( {\mathsf u}'{{\mathsf M}_A}^{\prime 1/2}{\mathsf M}_{B*}'-{\mathsf u}{\mathsf M}_A^{1/2}{\mathsf M}_{B*} ) \\
    &(|(\xi-\xi_*)\cdot n|) dn d\xi_* d\xi 
\end{align*}
which proves that $L_{AB}$ is self adjoint since the above expression is symmetric 
in ${\mathsf u} ,{\mathsf v}$. Furthermore putting $ {\mathsf u} ={\mathsf v} $ we get

$$ ( L_{AB}{\mathsf u} ,{\mathsf u} ) = -\frac{1}{2} \sigma^2_{AB}\int \frac{1}{{\mathsf M}_A{\mathsf M}_{B*}} {\left(
 {\mathsf u}'{{\mathsf M}_A'}^{1/2}{\mathsf M}_{B*}'-{\mathsf u} 
{\mathsf M}_A^{1/2}{\mathsf M}_{B*} \right)}^2 (|(\xi-\xi_*)\cdot n|) dn d\xi_* d\xi$$
which proves that 
$$ (L_{AB}{\mathsf u} ,{\mathsf u})  \leq 0 $$
as well as 
$$ (L_{AB} 
{\mathsf u} ,{\mathsf u} ) =0  \iff {\mathsf u}{\mathsf M}_A^{1/2}{\mathsf M}_{B*}
 ={\mathsf u}'{\mathsf M}_A^{\prime1/2}{\mathsf M}_{B*}' .$$
From ${\mathsf M}_A {\mathsf M}_{B*} = {\mathsf M}_A' {\mathsf M}_{B*}'$ and the above, one has that
 for all $\xi \in {\mathbb R}^3$
$$ {\mathsf u} /{\mathsf M}_A^{1/2} - {\mathsf u}'/{\mathsf M}_A^{\prime 1/2}=0. $$
As in the argument for the Boltzmann collision invarians if for all $\xi$ we have that $\phi(\xi)-\phi'(\xi)=0$ then $\phi$ is constant, which with the above yields that
$$ {\mathsf u} = c {\mathsf M}_A^{\frac{1}{2} } \text{ for } c \in {\mathbb R}. $$

Thus in contrast to the usual linearized Boltzmann operator $L_{AB}$ has a one
 dimensional kernel which corresponds to the conservation of mass law (obtained if one integrates
 the original equation \eqref{lseq1} against the constant 1). Notice that due to the dissipation 
of energy and momentum in a diffusion no further conservation laws could have been expected.

\end{proof}

Let us write $L_{AB}{\mathsf u}$ as follows:

$$L_{AB}{\mathsf u}  =  \frac{\sigma_{AB}^2}{{\mathsf M}_A^{1/2}}\int_{S^+}\int_{\mathbb{R}^3} 
{\mathsf u}'{{\mathsf M}_A'}^{1/2}{{\mathsf M}_B}'_*  |(\xi-\xi_*) \cdot n| d \xi_* dn  
       -\sigma_{AB}^2\int_{S^+}\int_{\mathbb{R}^3} {\mathsf u} {{\mathsf M}_B}_*  |(\xi-\xi_*) \cdot n| d \xi_* dn $$

    $$ = {\mathsf K}^{AB} {\mathsf u} - \nu^{AB} (|\xi|){\mathsf u} . $$
We have the following lemma summarizing the above and stating some properties of the
 ${\mathsf K}^{AB}-\nu^{AB}$ decomposition:

\begin{lemma}  \hfill \\
\label{lemma32}
\noindent (i) It immediately follows $\nu^{AB}= \frac{\sigma_{AB}^2  }{\sigma_{BB}^2}\nu$ where $\nu$ is exactly as in \eqref{In3}, so we have again

$$0 \leq c_1 \leq \nu^{AB} (|\xi|) \leq c_2 (1+|\xi|)$$
for some positive constants $c_1$ and $c_2$.

\noindent (ii) With the normalization $RT=1, \rho_B=1$ the operator $ {\mathsf K}^{AB} {\mathsf u} $ is given by 

$${\mathsf K}^{AB} {\mathsf u} (\xi) =  \int_{\mathbb{R}^3}k^{AB}(\xi,\xi_*) 
{\mathsf u} (\xi_*) d\xi_* \text{,   where} $$
$$k^{AB}(\xi,\xi_*)=  \frac{\sigma_{AB}^2 }{ \sqrt{8 \pi}    } \left( 1 
+ \frac{ m_A }{ m_B   } 
\right)^2
\frac{1}{|\xi_*-\xi|}
e^{-\frac{\|\xi_*|^2-|\xi|^2|^2}{8|\xi_*-\xi|^2} -\frac{m_A^2}{m_B^2}\frac{ |\xi_*-\xi|^2}{8} },$$
which is symmetric and uniformly square integrable in each variable.

\noindent (iii)The operator ${\mathsf K}^{AB}$ is compact $L^2(\mathbb{R}^3) \rightarrow L^2(\mathbb{R}^3)$. 

\noindent (iv) From i, ii, and iii we conclude that the operator $L_{AB}$ is a closed unbounded operator 
on $L^2(\mathbb{R}^3)$ with dense domain the set $D=\{ {\mathsf u} \text{ } |\text{ } (1+|\xi|){\mathsf u}
  \in L^2(\mathbb{R}^3)  \}$ and its kernel  is ${\mathsf E}_D \equiv{\mathsf M}_A^{1/2}$. 
Moreover $ L_{AB}$ is self-adjoint and non-positive, i.e.:

\begin{align}
\nonumber 
(L_{AB}{\mathsf u} ,{\mathsf v} )&=({\mathsf u} ,L_{AB} {\mathsf v}) \\
\nonumber 
(L_{AB} {\mathsf u},{\mathsf u} )& \leq 0
\end{align}
\noindent (v) It is strictly negative on functions orthogonal to the kernel, i.e. there exists a $\mu > 0$ such that if ${\mathsf u}  \perp {\mathsf E}_D$ then
$$ (L_{AB}{\mathsf u} ,{\mathsf u} ) \leq -\mu ({\mathsf u} ,{\mathsf u} ).$$
\end{lemma}

The proof is analogous to the proof for the full Boltzmann linearized operator. See \cite{Cerc}. \\ \\ 

For later use we will also need the following

\begin {lemma} \hfill \\
\label{lemmaba} 
\noindent (i) The operator $L_{BA}$ is a compact operator $L_{\xi}^2 \rightarrow L^2_{\xi}$ and is given as
$$L_{BA}f(\xi) =\frac{ 1 }{  \sqrt{{\mathsf M}_B}}  Q^{BA}( {\mathsf M}_B, {\mathsf f}_A \sqrt{{\mathsf M}_A} ) 
= \int k^{BA}(\xi,\xi_*) f(\xi_*) d\xi_*$$
where the kernel satisfies  $k^{BA}(-\xi,\xi_*)=k^{BA}(\xi,-\xi_*)$.

\noindent (ii)  The (unbounded) operator
$$\left(L_{AB}+ \frac{\sqrt{{\mathsf M}_B}}{  \sqrt{{\mathsf M}_A}  } L_{BA} \right) {\mathsf f} =
\frac{1}{  \sqrt{{\mathsf M}_A}  } Q^{AB}( \sqrt{{\mathsf M}_A} {\mathsf f},{\mathsf M}_B, )+
\frac{1}{  \sqrt{{\mathsf M}_A}  } Q^{BA}({\mathsf M}_B, \sqrt{{\mathsf M}_A}  {\mathsf f} ) $$
is orthogonal to all collision invariants $\sqrt{{\mathsf M}_A}, \xi_i\sqrt{{\mathsf M}_A}, |\xi|^2\sqrt{{\mathsf M}_A}$. 
\end{lemma}

\begin{proof}

(i) The kernel $k^{BA}$ is obtained analogously to that in the previous lemma.
(ii) To see this conservative property denote by $'AB$ and $'BA$ the velocity transformations appearing in the two collisional integrals. Then more generally observe that
$$\langle Q^{AB}(f,h)(\xi)+Q^{BA}(h,f)(\xi), \phi(\xi) \rangle =$$
$$=\sigma_{AB}^2  \int_{\mathbb{R}^3} \int_{\mathbb{R}^3} \int_{S^+}
\phi(f'^{AB}h'^{AB}_* - fh_* )|(\xi-\xi_*) \cdot n| + \phi(h'^{BA}f'^{BA}_* - hf_* ) |(\xi-\xi_*) \cdot n| dn d \xi_*  d\xi= $$
$$=\sigma_{AB}^2  \int_{\mathbb{R}^3}\int_{\mathbb{R}^3} \int_{S^+} (L + R)\phi dn d \xi_* d\xi = \frac{\sigma_{AB}^2} {2} \int_{\mathbb{R}^3}\int_{\mathbb{R}^3} \int_{S^+} (L\phi + R\phi+L\phi+R\phi) dn d \xi_* d\xi. $$

Now notice that the transformation $T_*: \xi \rightarrow \xi_*, \xi_*\rightarrow \xi$ takes $R$ to $L$ and reversely. Moreover the transformation $T'^{AB}:\xi \rightarrow \xi'^{AB}, \xi_* \rightarrow {\xi'}_*^{AB}$ takes $L$ to $- L$. So after performing $T_*$ on the second term, $T'^{AB}$ on the third term and $T'^{AB}T_*$ on the fourth term we get

$$\langle Q^{AB}(f,h)(\xi)+Q^{BA}(h,f)(\xi), \phi(\xi) \rangle = $$
$$\frac{\sigma_{AB}^2} {2} \int_{\mathbb{R}^3}\int_{\mathbb{R}^3}\int_{S^+} L(\phi -\phi'^{AB}+\phi_*-\phi_*'^{AB}) dn d \xi_* d\xi.$$
From this as in the treatment of the usual collision operator it follows the above expression is 0 if $\phi(\xi)$ is a linear combination of the functions $1,\xi_i, |\xi|^2$ and the Lemma follows.
\end{proof}

\noindent{\bf Macro-Micro Decomposition for Boltzmann Diffusive equation } \\
Similar to the macro-micro decomposition $({\mathsf P}_0,{\mathsf P}_1)$ for $L$, 
from Lemma \ref{lemma31} and Lemma \ref{lemma32}, we can decompose $L^2_\xi$ into 
$ker(L_{AB}) \oplus ker(L_{AB})^\perp$ as follows
\begin{equation}
\label{mmd}
\begin{cases}
{\mathsf u} = {\mathsf P}_0^D {\mathsf u} + {\mathsf P}_1^D {\mathsf u}, \\
{\mathsf P}_0^D  \equiv {\mathsf E}_D  \otimes \left \langle  {\mathsf E}_D 
\right|,
\; \;
( {\mathsf P}_0^D {\mathsf u}  \equiv ({\mathsf u}, {\mathsf E}_D  ) {\mathsf E}_D ) \\
{\mathsf P}_1^D \equiv 1- {\mathsf P}_0^D,\\
{\mathsf E}_D \equiv {\mathsf M}_A^{1/2}, \; \;( ker(L_{AB}) = span\{{\mathsf E}_D\}).
\end{cases} 
\end{equation} 
\end{subsection}
\begin{subsection}{The Spectrum Decomposition of Boltzmann Diffusive equation } \hfill \\

Taking the Fourier Transform in the variable $x$ of the equation \eqref{lseq1} we obtain the equation 
$$
\partial_t \hat {\mathsf g} = ( L_{AB} - i k \xi^1)  \hat {\mathsf g} 
= {\mathsf K}^{AB} \hat {\mathsf g}   -\nu \hat  {\mathsf g}  - i k \xi^1 \hat {\mathsf g}
$$ 
which leads us to study the spectral properties of the operator $L_{AB}-i k\xi^1$.

Similarly to the explanation in \cite{Cerc} we have the following lemma.
\begin{lemma}

For each fixed $k \in \mathbb{R}$, the spectrum  of $L_{AB}-ik\xi^1$ in the domain 
 $Re \lambda > - \nu_0$  consists of isolated eigenvalues with non-positive real parts. 

\end{lemma}

Next we want to consider the behavior of the spectrum for $|k| \ll 1$. We have the following lemma:

\begin{lemma}
\hfill \\ 
 \label{spectheo} 
\noindent (i) Given $\kappa_0>0$, if $|k|>\kappa_0>0$, $k \in \mathbb{R}$ then there exists $ \epsilon(\kappa_0) > 0$ such that if $\lambda$ 
is an eigenvalue of $L_{AB}-ik\xi^1$ then $Re \lambda < - \epsilon$.  \\
\noindent (ii) There exist $\kappa_0$ and $\delta$ such that for $|k|<\kappa_0$ the spectrum with
 $|\lambda(k)|<\delta$ is a single point which converges to the origin and it , together with
 its eigenvector ${\mathsf e}_D(k)$ is an analytic function of $k$ for $|k|<\kappa_0$.  \\
\noindent (iii) The expansions of the eigenvector and eigenvalue are as follows:
\begin{align}
\label{spin1}
\lambda(k) &= -a_2k^2+ 0(k^3) \text{  where $a_2$ is a positive real number, and } \\
{\mathsf e}_D(k)&={\mathsf E}_D+ {\mathsf e}_{D}'(0)k+O(1)|k|^2, \text{ where }
 {\mathsf e}_D'(0)= iL_{AB}^{-1}\xi^1 {\mathsf E}_D.
\label{spin2}
\end{align}
\end{lemma} 

\begin{proof}
Statements (i) and (ii) follow from standard perturbation theory and the fact that $\nu^{AB}(|\xi|)<c(1+|\xi|)$, 
see \cite{Cerc}.

For each fixed $k$, the eigenvalue problem for the operator $-i k \xi^1 + L_{AB}$ is 
\begin{equation}
\label{eigen_a}
(-i k \xi^1 +L_{AB}) {\mathsf e}_D(k) = \lambda(k) {\mathsf e}_D(k).
\end{equation} 
Now apply the Macro-Micro decomposition $({\mathsf P}_0^D, {\mathsf P}_1^D)$ given in \eqref{mmd} to both 
\eqref{eigen_a} and ${\mathsf e}_D(k)$ 
 to get  \begin{equation}
\label{eigen_b}
\begin{cases}
\displaystyle  
{\mathsf e}_D(k) = {\mathsf P}_0^D {\mathsf e}_D(k) +
  {\mathsf P}_1^D {\mathsf e}_D(k) \equiv \psi_0 + \psi_1, \\
ik {\mathsf P}^D_0\xi^1(\psi_0+ \psi_1)=\lambda\psi_0, \\
L_{AB} \psi_1 + ik {\mathsf P}^D_1\xi^1(\psi_0+ \psi_1)=\lambda\psi_1.
\end{cases}  
 \end{equation}

The last equation allows us to express
\begin{equation} \label{psi1} \psi_1=ik[L_{AB}-ik {\mathsf P}^D_1\xi^1-\lambda]^{-1}\psi_0. \end{equation}
Which substituted into the second one in \eqref{eigen_b} with the property 
$({\mathsf E}_D, \xi^1 {\mathsf E}_D)=0$ gives
$$k^2 {\mathsf P}^D_0\xi^1[L_{AB}-ik {\mathsf P}_1\xi^1-\lambda]^{-1}\psi_0 = \lambda\psi_0.$$
Taking inner product with ${\mathsf M}_A^{1/2}$ we get the following implicit relation for $\lambda(k)$
$$\Delta(\lambda,k)= \int {\mathsf E}_D( k^2\xi^1[L_{AB}-ik {\mathsf P}^D_1\xi^1-\lambda]^{-1}-\lambda) 
{\mathsf E}_Dd\xi=0$$
We note that 
$$ \left. {{\partial_k \Delta}}\right\vert_{\lambda,k=0}=0, \;  \;  \;
\left. {{\partial_\lambda  \Delta}{}} \right\vert_{\lambda,k=0}=1 $$
$$\left. {{\partial^2_k \Delta}{}}\right|_{\lambda,k=0}= \int \xi^1 {\mathsf E}_DL_{AB}^{-1}\xi^1
{\mathsf E}_D d \xi =-a_2 <0$$
Therefore by the implicit function theorem we get 
\begin{align}
\nonumber 
\left. { \frac{d \lambda}{ d k}} \right\vert_{ k=0}   =  \left.
   -\frac {{  {\partial_k \Delta}{}}}
    { { {\partial_\lambda  \Delta}{}} } \right\vert_{\lambda=0, k=0}=0.
 \\
\nonumber 
 \left.
{\frac {d^2 \lambda}{dk^2}} \right|_{k=0}   = 
   -  \left. \frac {{ {\partial^2_k \Delta}{ }}}
    { { {\partial^2_\lambda  \Delta}{ }}} \right\vert_{\lambda=0, k=0} = a_2>0 .
\end{align}
Next we consider the eigenvector ${\mathsf e}^D(k)$ near $k=0$. Normalizing
 ${\mathsf e}^D(k)={\mathsf E}_D+\psi_1(k)$ from the expression \eqref{psi1} after
 differentiation we immediately obtain the desired expansion.
\end{proof}

Note that the number $a_2 = - (\xi^1 {\mathsf E}_D , L^{-1}_{AB} \xi^1 {\mathsf E}_D) $ is realized as {\bf  cross species diffusion coefficient}.

\end{subsection}

\end{section}

\begin{section}{On the Linear Solution Operators for the Boltzmann Diffusive equation }
Consider the initial value problem
\begin{equation}
\label{bde}
\begin{cases}
\partial_t {\mathsf g} + \xi^1 \partial_x {\mathsf g} = L_{AB}{\mathsf g},\\
{\mathsf g}\vert_{t=0} = {\mathsf g}_{in},
\end{cases} 
\end{equation} 
where ${\mathsf g}_{in}$ satisfies \eqref{indata}. 
Then, we have
$$ \hat {\mathsf g}(k,t)= e^{(-i\xi^1 k + L_{AB})t} \hat {\mathsf g}_{in}(k).$$
 Now taking inverse Fourier transform as before we have for the Green's function
$$ \mathbb{G}_{AB}(x,t) \equiv \frac {1}{2\pi} \int_{\mathbb{R}}e^{ikx+(-i\xi^1 k + L_{AB})t} dk.$$

Similarly to the Long Wave-Short Wave decomposition \eqref{l-s-de}, we also 
use the same decomposition for ${\mathbb G}_{AB}$:

$$ 
\begin{cases}
\mathbb{G}_{AB}(x,t) = \mathbb{G}_{AB;L}(x,t)+\mathbb{G}_{AB;S}(x,t),  \\
\displaystyle   \mathbb{G}_{AB;L}(x,t) 
\equiv  \frac {1}{\sqrt{2\pi}} \int_{|k|<\kappa_0}e^{ikx+(-i\xi^1 k + L_{AB})t} dk,\\
\displaystyle  
\mathbb{G}_{AB;S}(x,t) \equiv 
\frac {1}{\sqrt{2\pi}} \int_{|k|>\kappa_0}e^{ikx+(-i\xi^1 k + L_{AB})t} dk.
\end{cases}
$$
Furthermore, in the case of $\mathbb{G}_{AB}$, if $\Pi_k^D$ is the projection on the eigenvector
 ${\mathsf e}_D(k)$ corresponding to the eigenvalue $\lambda(k)$ discussed in Lemma \ref{spectheo} 
we can write
$$ \mathbb{G}_{AB;L}(x,t) = \frac {1}{\sqrt{2\pi}} \int_{|k|<\kappa_0}e^{ikx+(-i\xi^1 k + L_{AB})t} (\Pi_k^D+\Pi_k^{D\perp})dk 
 \equiv
 \mathbb{G}_{AB;L,0}(x,t)+ \mathbb {G}_{AB;L,\perp}(x,t).$$

Similar to \eqref{spec_d2}, the operator $\Pi_k^{D\perp}$ satisfies that for 
$|k| \le \kappa_0$,
\begin{equation}
\label{spec_d3}
\|  e^{ (- i k \xi^1 +L_{AB})t}\Pi_k^{D\perp} \|_{L^2_\xi} \le C e^{-t/C} \text{ for some }C>0.
\end{equation} 
Writing ${\mathsf e}_D(k)={\mathsf E}_D + {\mathsf e}_D'(0)k + O(1)k^2$, we have that
$$\Pi_{k}^D = {\mathsf E}_D \otimes \left \langle {\mathsf E}_D \right| + k {\mathsf e}_D'(0) \otimes 
\left \langle  {\mathsf E}_D \right| + k 
{\mathsf E}_D \otimes \left \langle {\mathsf e}_D'(0) \right| + O(1) k^2. $$

>From property (iii) of Lemma \ref{spectheo}, one has that
\begin{equation}
\label{spforg}
e^{(-ik \xi^1 +L)t} \Pi_k^D 
= e^{-\lambda(k) t}
\left(
 {\mathsf E}_D
\otimes \left \langle {\mathsf E}_D
\right| + k 
( {\mathsf e}_D'(0)
\otimes \left \langle {\mathsf E}_D
\right| 
+  {\mathsf E}_D
\otimes \left \langle {\mathsf e}_D'(0)
\right|
 ) +O(k^2)\right).
\end{equation} 
Here in \eqref{spforg}, all functions are analytic around
$k=0$. \\
\begin {theorem} \label{basictheo}
Recall that we denote $|||.|||:=\|.\|_{L^{\infty}_x(L^\infty_{\xi,3})}$ and 
$({\mathsf P}_0^D, {\mathsf P}_1^D)$ the macro-micro decomposition given
in \eqref{mmd}.
 The Boltzmann Diffusion Green function satisfies the estimates:

In the region $|x|<2|\lambda_1|t$,
\begin{align} \label{basictheoa}
\|\mathbb{G}_{AB}{\mathsf g}_{in} (x,t)\|_{L^2_{\xi}}& = O(1)|||{\mathsf g}_{in}|||
\left( \frac{ e^{ - \frac{ (x-\lambda_2 t)^2 }{  C(1+t)  } } }{ \sqrt{1+t}   }  +e^{-(t+|x|)/C}\right).\\
 \label{basictheob}
\| {\mathsf P}_1^D\mathbb{G}_{AB} {\mathsf g}_{in} (x,t)\|_{L^2_{\xi}} \text {,  } \|\mathbb{G}_{AB}
 {\mathsf P}_1^D {\mathsf g}_{in} (x,t)\|_{L^2_{\xi}} &= 
 O(1)|||{\mathsf g}_{in}|||
\left( 
 \frac{ e^{ - \frac{ (x-\lambda_2 t)^2 }{  C(1+t)  } } }{ (1+t)   }   + e^{-(t+|x|)/C} \right),\\ \label{basictheoc}
\|{\mathsf P}_1^D \mathbb{G}_{AB} {\mathsf P}_1^D {\mathsf g}_{in} 
(x,t) \|_{L^2_{\xi}} &= 
 O(1)|||{\mathsf g}_{in}|||
\left( \frac{ e^{ - \frac{ (x-\lambda_2 t)^2 }{  C(1+t)  } } }{ (1+t)^{3/2}   }  +e^{-(t+|x|)/C}\right).
\end{align}

In the region $|x|\geq 2|\lambda_1|t$:
$$\|\mathbb{G}_{AB}{\mathsf g}_{in} (x,t)\|_{L^2_{\xi}} = O(1)|||{\mathsf g}_{in}|||e^{-(t+|x|)/C}.$$

\end{theorem}
\begin{proof}
First, we claim the Particle-Wave decomposition in \eqref{p-w-d} is 
also valid for  \eqref{bde}, since the operator ${\mathsf K}^{AB}$ for
$L_{AB}$  and 
${\mathsf K}$ for $L$ both have the properties necessary for the proofs in \cite{liuyug1}  for 
obtaining Lemma \ref{lemma-p-w-d}. Thus we have that 
there exist $({\mathbb P}_j^D, {\mathbb W}_j^D)$ satisfying
\begin{equation}
\label{p-w-d-2}
 \begin{cases}
{\mathsf g}(x,t) = {\mathbb P}_j^D(x,t) + {\mathbb W}_j^D(x,t) \text{ for all }
x \in {\mathbb R}, t\ge 0, \\
{\mathbb P}_j^D(x,0) \equiv {\mathsf g}_{in}(x), \\
\|{\mathbb P}_j^D(x,t) \|_{L^2_\xi} \le C_j e^{-(|x|+t)/C_j} \text{ for all } x 
\in {\mathbb R}, t \ge 0, \\
\| (\partial_t + \xi^1 \partial_x - L_{AB}){\mathbb P}_j^D \|_{H^j_x(L^2_\xi)} 
\le C_j e^{-t/C_j}, 
\\
\|{\mathbb W}_j^D\|_{H^j_x(L^2_\xi)} \le C_j \text{ for all } t \ge 0 \text{(uniformly bounded in time)}.
\end{cases} 
\end{equation}
Here, we can choose $j=3$. 

>From property (i) of Lemma \ref{spectheo}, one can have 
that 
\begin{equation}
\label{gse}
\| {\mathbb G}_{AB;S}^t {\mathsf g}_{in}\|_{L^2_x(L^2_\xi)} \le
O(1) \|| {\mathsf g}_{in} \|| e^{-t/C} \text{ for some }C>0.
\end{equation} 
>From property (iii) of Lemma \ref{spectheo}, one can have
that 
\begin{equation}
\label{gle}
\| {\mathbb G}_{AB;L}^t {\mathsf g}_{in}\|_{H^3_x(L^2_\xi)} \le
O(1) \|| {\mathsf g}_{in} \||.
\end{equation} 
By combining \eqref{p-w-d-2} and \eqref{gle}, one has that
\begin{equation}
\label{gse1}
\| {\mathbb G}_{AB;S} {\mathsf g}_{in} - {\mathbb P}_j^D \|_{H^3_x(L^2_\xi)} =
\| {\mathbb G}_{AB;L} {\mathsf g}_{in} - {\mathbb W}_j^D \|_{H^3_x(L^2_\xi)} =
  O(1) \|| {\mathsf g}_{in} \|| \text{ for all } t\ge 0.
\end{equation} 
>From \eqref{p-w-d-2}, \eqref{gse},
\eqref{gse1}, and Sobolev's inequality,
one has that
\begin{equation}
\label{diff}
\| {\mathbb G}_{AB;S} {\mathsf g}_{in}  \|_{L^\infty_x(L^2_\xi)} \le O(1) e^{-t/C} \text{ for some } C>0.
\end{equation} 
Then, by the spectral property \eqref{spforg} and \eqref{spin1}, one can apply 
a complex contour integral similar to that
in \cite{liuyug1} to obtain that
\begin{equation}
\label{diff_2}
\| {\mathbb G}_{AB;L}{\mathsf g}_{in}(x,t)
\|_{L^2_\xi} \le O(1) \frac{ 
e^{- \frac{ (x- \lambda_2 t)^2 }{C (t+1)    }  }   }{  \sqrt{(1+t)}  }
\text{ for } |x| \le 2 |\lambda_1 |t.
\end{equation} 
Hence, \eqref{diff} and \eqref{diff_2}
conclude \eqref{basictheoa}. The estimates
\eqref{basictheob} and \eqref{basictheoc}
will follow by taking expansion of higher
order terms in $k$  in 
\eqref{spforg} into account.

In the region $|x| \ge 2 |\lambda_1| t$, 
one can apply the  
weighted energy estimates 
 in \cite{liuyug1}
to the wave component ${\mathbb W}^D_3(x,t)$.
Then, one can show that 
$\|{\mathbb W}^D_3(x,t)\|_{L^2_\xi}$ decays 
to zero exponentially in $x$ and $t$
in the region $|x| \ge 3 |\lambda_1| t/2$.
Then, this and \eqref{p-w-d-2} will conclude this theorem for $|x| > 2 | \lambda_1 | t$. 

\end{proof} 

\begin{remark}
Since the decomposition  $({\mathbb P}_j^D, {\mathbb W}_j^D)$  satisfies
\eqref{p-w-d-2}  one can show that for large times it is equivalent to the decomposition $({\mathbb G}_{AB;L}^t, {\mathbb G}_{AB;S}^t)$ due to \eqref{diff}. We call this procedure {\bf ``separation of scales''}.
\end{remark} 
\end{section}

\begin{section}{On Cross Species Mass Diffusion and 
Sound Wave Interaction } 
Consider the initial value problem \eqref{upper} with initial data given
in \eqref{indata}:
\begin{equation}
 \label{lynsis} \begin{array}{l}
\begin{cases} 
\partial_t {\mathsf g} + \xi^1 \partial_x {\mathsf g}  = L_{AB} {\mathsf g}, \\ 
\partial_t {\mathsf h}  + \xi^1 \partial_x {\mathsf h}  = L{\mathsf h} +  L_{BA}{\mathsf g},
\end{cases} \\ \; \;
({\mathsf g},{\mathsf h})\vert_{t=0} = ({\mathsf g}_{in}, {\mathsf h}_{in}).
\end{array}
\end{equation}
The solution of the first equation is $\mathbb{G}_{AB}^t {\mathsf g}_{in}$ and so satisfies the estimate provided by Theorem \ref{basictheo}. 
Next, we continue to study the behavior of ${\mathsf h}(x,t)$ just treating
$L_{BA} {\mathsf g}$ as a source term to the equation for ${\mathsf h}$. 
However, the function ${\mathsf g}$ and 
the operator $\partial_t + \xi^1 \partial_x - L $  have a resonance around the
trajectory $x= \lambda_2 t=0$. This makes this system very interesting. 

\noindent {\bf 1)  Estimate in the Finite Mach Region $x<|2\lambda_1t|$}

By Duhamel's principle one can represent the 
solution ${\mathsf h}(x,t)$ as follows
\begin{multline}
\nonumber 
{\mathsf h}(x,t) = \int_{{\mathbb R}} {\mathbb G}(x-y,t) {\mathsf h}_{in}(y) dy \\ +
\left(  \int_0^{t^\frac{1}{2} }+\int_{t^{\frac{1}{2}}}^{t-t^\frac{1}{2} } + \int_{t-t^\frac{1}{2}}^t
\right) \int_{{\mathbb R}^2} 
{\mathbb G}(x-y,t-s)   
L_{BA} {\mathbb G}_{AB}(y-z,s) {\mathsf g}_{in}(z)
dz dy ds.
\end{multline}

Due to  \eqref{ABBA} and the above, we have the following {\bf ``microscopic cancellation''}
\begin{equation}
\label{51b}
 L_{BA}
 {\mathsf J} =
\frac{1}{  \sqrt{{\mathsf M}_B}  } Q^{BA}({\mathsf M}_B,
\sqrt{{\mathsf M}_A} {\mathsf P}_1^D
{\mathsf J} ) 
=  L_{BA} {\mathsf P}_1^D
{\mathsf J}.
\end{equation}
With this microscopic cancellation, \eqref{basictheob}, and Theorem \ref{mainA}, one has that
\begin{multline}
\nonumber 
\left\| 
 \int_{t-t^\frac{1}{2}}^t
\int_{{\mathbb R}^2} 
{\mathbb G}(x-y,t-s)   
L_{BA} {\mathbb G}_{AB}(y-z,s) {\mathsf g}_{in}(z)
dz dy ds
\right \|_{L^2_\xi} \\ =
\left\| 
 \int_{t-t^\frac{1}{2}}^t
\int_{{\mathbb R}^2} 
{\mathbb G}(x-y,t-s)   
L_{BA} {\mathsf P}_1^D {\mathbb G}_{AB}(y-z,s) {\mathsf g}_{in}(z)
dz dy ds 
\right \|_{L^2_\xi}  \\ \le 
O(1) \sum_{j=1}^3 \int_{t-t^\frac{1}{2}}^t
\int_{{\mathbb R}} 
\left( \frac{  e^{ - \frac{ (x-y - \lambda_j(t-s))^2 }{  C(t-s)  } }  }{  \sqrt{1+(t-s)}  } 
+ e^{-(|x-y|+(t-s))/C} \right)
\left( \frac{  e^{ - \frac{ y^2 }{  Cs  } }  }{  (1+s)  } 
+ e^{-(|y|+s)/C} \right)
dy  ds \\ =
O(1) \sum_{j=1}^3
\left( \frac{  e^{ - \frac{ (x - \lambda_jt)^2 }{  Ct  } }  }{  \sqrt{(1+t)}  } 
+ e^{-(|x|+t)/C} \right).
\end{multline} 
Furthermore with the microscopic cancellation
\eqref{51b},  one can rewrite 
$ L_{BA} {\mathsf P}_1^D {\mathsf g}$ in \eqref{lynsis} through 
$\partial_t {\mathsf g} + \xi^1 \partial_x {\mathsf g} 
- L_{AB} {\mathsf g} =0$:
\begin{equation}
\label{51c}
 L_{BA} {\mathsf P}_1^D {\mathsf g} = 
 L_{BA}
L_{AB}^{-1}  {\mathsf P}_1^D(\partial_t + \xi^1 \partial_x)
{\mathsf g}.
\end{equation} 
This gives
\begin{multline}
\label{51d}
\left\| 
 \int^{t^\frac{1}{2}}_0
\int_{{\mathbb R}^2} 
{\mathbb G}(x-y,t-s)   
L_{BA} {\mathbb G}_{AB}(y-z,s) {\mathsf g}_{in}(z)
dz dy ds
\right \|_{L^2_\xi} \\ =
\left\| 
 \int^{t^\frac{1}{2}}_0
\int_{{\mathbb R}^2} 
{\mathbb G}(x-y,t-s)   
L_{BA} L_{AB}^{-1} {\mathsf P}_1^D 
(\partial_s + \xi^1 \partial_y)
{\mathbb G}_{AB}(y-z,s) {\mathsf g}_{in}(z)
dz dy ds 
\right \|_{L^2_\xi}  
\\
  \le 
\left\| 
 \int^{t^\frac{1}{2}}_0
\int_{{\mathbb R}^2} 
\partial_s 
{\mathbb G}(x-y,t-s)   
L_{BA} L_{AB}^{-1} {\mathsf P}_1^D 
{\mathbb G}_{AB}(y-z,s) {\mathsf g}_{in}(z)
dz dy ds 
\right \|_{L^2_\xi} 
\\
+
\left\| 
 \int^{t^\frac{1}{2}}_0
\int_{{\mathbb R}^2}  \partial_y
{\mathbb G}(x-y,t-s)   
L_{BA}  L_{AB}^{-1} {\mathsf P}_1^D 
 \xi^1 
{\mathbb G}_{AB}(y-z,s) {\mathsf g}_{in}(z)
dz dy ds 
\right \|_{L^2_\xi}
\\ +
\left\| 
 \left.
\int_{{\mathbb R}^2}  
{\mathbb G}(x-y,t-s)   
L_{BA}  L_{AB}^{-1} {\mathsf P}_1^D 
{\mathbb G}_{AB}(y-z,s) {\mathsf g}_{in}(z)
dz dy  \right\vert_{s=0}^{s=t^{\frac{1}{2}}}
\right \|_{L^2_\xi}
  \\ \le 
O(1) \sum_{j=1}^3 \int^{t^\frac{1}{2}}_0
\int_{{\mathbb R}} 
\left( \frac{  e^{ - \frac{ (x-y - \lambda_j(t-s))^2 }{  C(t-s)  } }  }{  {1+(t-s)}  } 
+ e^{-(|x-y|+(t-s))/C} \right)
\left( \frac{  e^{ - \frac{ y^2 }{  Cs  } }  }{  \sqrt{(1+s)}  } 
+ e^{-(|y|+s)/C} \right)
dy  ds \\
+ 
O(1) \sum_{j=1}^3
\left( \frac{  e^{ - \frac{ (x - \lambda_jt)^2 }{  Ct  } }  }{  \sqrt{(1+t)}  } 
+ e^{-(|x|+t)/C} \right) 
  =
O(1) \sum_{j=1}^3
\left( \frac{  e^{ - \frac{ (x - \lambda_jt)^2 }{  Ct  } }  }{  \sqrt{(1+t)}  } 
+ e^{-(|x|+t)/C} \right),
\end{multline} 
where we noted that $L_{BA}  L_{AB}^{-1} {\mathsf P}_1^D$ is a bounded operator on $L^2_{\xi}$ since $L_{AB}^{-1}{\mathsf P}_1^D$ is bounded, and $L_{BA}$ is compact. Also as is clear by the proof of Theorem \ref{mainA} we can justify the differentiation of $\mathbb{G}$ since in the finite Mach region $\mathbb{G}$ equals $\mathbb{G}_L$ up to a term with exponential decay in time and space. Now in the finite Mach region $\mathbb{G}_L$ is a smooth function and the estimates on the derivatives follow from the explicit formula used in the proof of Theorem \ref{mainA}.

Denote the remaining integral by $ {\mathsf q}(x,t)$: 
\begin{multline} 
 {\mathsf q}(x,t) g_{in}\equiv
\int^{t-t^\frac{1}{2}}_{t^{\frac{1}{2} }}
\int_{{\mathbb R}^2} 
{\mathbb G}(x-y,t-s)   
L_{BA}  {\mathsf P}_1^D{\mathbb G}_{AB}(y-z,s) {\mathsf g}_{in}(z)
dz dy ds =\\  
 =\int_{t^{\frac{1}{2}}}^{t-  t^{\frac{1}{2}} } \int_{-\infty}^{\infty} e^{i k x }  e^{(-i k \xi^1 +L )(t-s)}   L_{BA}
{\mathsf P}_1^D
e^{(-i k \xi^1 +L_{AB} )s} {\mathsf g}_{in}(k) dk ds = \int^{t-t^\frac{1}{2}}_{t^{\frac{1}{2} }} \int_{-\infty}^{\infty} e^{i k x } \hat {\mathsf q}(k, t)dkds
\end{multline}
For this function ${\mathsf q}(x,t)$, we consider the two decompositions:
\begin{equation}
\label{two-dec}
\begin{cases}
{\mathsf q}(x,t) \equiv {\mathsf q}_L(x,t) + {\mathsf q}_S(x,t), \text{ Long Wave-Short Wave Decomposition}, \\
{\mathsf q}(x,t) \equiv {\mathsf q}_{\mathscr P}(x,t) + {\mathsf q}_{\mathscr W}(x,t), 
\text{ Particle-Wave Decomposition},
\end{cases} 
\end{equation} 
where 
\begin{equation}
\label{two-dec_d}
\begin{cases}
\displaystyle  {\mathsf q}_L(x,t) \equiv \int_{\sqrt{t}}^{t- \sqrt{t}} \int_{|k| \le \kappa_0} e^{i k x } \hat {\mathsf q}(k,t) dx, \; \;\\ 
{\mathsf q}_S(x,t) \equiv {\mathsf q}(x,t) - {\mathsf q}_L(x,t), \\
\displaystyle  
{\mathsf q}_{{\mathscr P}}(x,t) \equiv \int_{\sqrt{t}}^{t- \sqrt{t}} \int_{\mathbb R}{\mathbb P}_j(x-y,t-s) 
 L_{BA} {\mathsf P}_1^D {\mathbb P}_j^D(y,s) dy
 ds, \; \;\\
{\mathsf q}_{\mathscr W}(x,t)  \equiv  {\mathsf q}(x,t)  - {\mathsf q}_{\mathscr P}(x,t) =\\ 
\int^{t-t^\frac{1}{2}}_{t^{\frac{1}{2} }} {\mathbb P}_j(t-s)\ast  L_{BA} {\mathsf P}_1^D{\mathbb W}_j^D(s)+ {\mathbb W}_j(t-s)\ast  L_{BA} {\mathsf P}_1^D{\mathbb P}_j^D(s)+ {\mathbb W}_j(t-s)\ast  L_{BA} {\mathsf P}_1^D {\mathbb W}_j^D ds
\end{cases} 
\end{equation} 
where ${\mathbb  P}_j$,  ${\mathbb  W}_j$ and ${\mathbb P}_j^D$, ${\mathbb  W}_j^D$ are given in \eqref{re-on-1} and \eqref{p-w-d-2} and the convolutions are in the space variable. \\ \\
The  consequence of these two decompositions is that  for some $C>0,$
\begin{equation}
\label{cons-tw}
\begin{cases}
\|{\mathsf q}_{\mathscr P}(x,t)\|_{L_\xi^2} \le C e^{-(|x|+t)/C} \text {    as a convolution of decaying exponentials in time and space}, \\
\| {\mathsf q}_{\mathscr W}\|_{H^3_x(L^2_\xi)} \le C (1+t) \text {    as a convolution with a regular function bounded in time} , \\
\| {\mathsf q}_L\|_{L^2_x(L^2_\xi)} \le C (1+t) \text {    due to the spectral properties }, \\
\| {\mathsf q}_S\|_{L^2_x(L^2_\xi)} \le C e^{-(1+t)/C} \text {    due to the spectral properties}.
\end{cases} 
\end{equation} 
Combining the representations in \ref{two-dec} and using \ref{cons-tw}:
\begin{equation}
\nonumber 
\begin{cases}
\| {\mathsf q}_S - {\mathsf q}_{\mathscr P} \|_{H^3_x(L^2_\xi)} = \| {\mathsf q}_{\mathscr W} - {\mathsf q}_L\|_{H^3_x(L^2_\xi)}  \le
C (1+t), \\ 
\| {\mathsf q}_S - {\mathsf q}_{\mathscr P} \|_{L^2_x(L^2_\xi)} \le 
C e^{-t/C}
\end{cases}
\end{equation}
Hence, this and Sobolev's inequality give
\begin{equation}
\nonumber 
\| {\mathsf q}_S \|_{L^\infty_x(L^2_\xi)} \le
O(1) e^{-(1+t)/C} \text{ for some } C>0.
\end{equation} 
which in the finite Mach region $x<|2\lambda_1t|$ gives the desired
\begin{equation}
\label{cons-tw2}
\| {\mathsf q}_S(x,t) \|_{L^2_\xi} \le
O(1) e^{-(t+|x|)/C} \text{ for some } C>0.
\end{equation} 

For $k<\kappa_0$ we consider a spectral decomposition of 
$ e^{(-i k \xi^1 +L)(t-s)} L_{BA} {\mathsf P}_1^D
e^{(-i k \xi^1 +L_{AB} )s} $ in terms of the spectral decompositions of the 
operators $-i k \xi^1 +L$  in \eqref{spec_d1} and $-i k \xi^1 +L_{AB}$ in \eqref{spforg}
as follows
\begin{multline}
\nonumber 
e^{(-i k \xi^1 +L)(t-s)} L_{BA} {\mathsf P}_1^D
e^{(-i k \xi^1 +L_{AB} )s} \\  = 
\sum_{j=1}^3 
e^{\sigma_j(k)(t-s) +\lambda(k) s}
 {\mathsf e}_j(k) \otimes \left \langle {\mathsf e}_j(k) \right|
L_{BA} {\mathsf P}_1^D
 {\mathsf e}_D(k) \otimes \left \langle {\mathsf e}_D(k) \right|
\\ +
\sum_{j=1}^3 
e^{\sigma_j(k)(t-s)}
 {\mathsf e}_j(k) \otimes \left \langle {\mathsf e}_j(k) \right| 
L_{BA} {\mathsf P}_1^D e^{(-i k \xi^1 +L_{AB})s} \Pi_k^{D\perp} \\
+
e^{\lambda(k)s } e^{(-i k \xi^1 +L)(t-s)} \Pi_k^{\perp}
L_{BA} {\mathsf P}_1^D 
 {\mathsf e}_D(k) \otimes \left \langle {\mathsf e}_D(k) \right| \\
+
 e^{(-i k \xi^1 +L)(t-s)} \Pi_k^{\perp}
L_{BA} {\mathsf P}_1^D e^{(-i k \xi^1 +L_{AB})s} \Pi_k^{D\perp}.
\end{multline}
>From \eqref{normal_} and \eqref{spin2}, it follows
\begin{multline} \label{52}
 {\mathsf e}_j(k) \otimes \left \langle {\mathsf e}_j(k) \right|
L_{BA} {\mathsf P}_1^D
 {\mathsf e}_D(k) \otimes \left \langle {\mathsf e}_D(k) \right|=
k {\mathsf E}_j \otimes \left \langle {\mathsf E}_j \right|
L_{BA} {\mathsf P}_1^D
 {\mathsf e}_D'(0) \otimes \left \langle {\mathsf E}_D \right| +O(1)k^2 \\
= -ik 
 {\mathsf E}_j \otimes \left \langle {\mathsf E}_j \right|
L_{BA} 
 L_{AB}^{-1} {\mathsf P}_1^D \xi^1 {\mathsf E}_D \otimes \left \langle {\mathsf E}_D \right| +O(1) k^2.
\end{multline} 

\begin{remark}
The components
$
e^{\sigma_j(k)(t-s) +\lambda(k) s}
 {\mathsf e}_j(k) \otimes \left \langle {\mathsf e}_j(k) \right|
L_{BA} {\mathsf P}_1^D
 {\mathsf e}_D(k) \otimes \left \langle {\mathsf e}_D(k) \right|
$ for $j=1,3$ are the Fourier transformation of the {\bf 
interaction of sound wave and cross species mass diffusion}.

The component $
e^{\sigma_2(k)(t-s) +\lambda(k) s}
 {\mathsf e}_2(k) \otimes \left \langle {\mathsf e}_2(k) \right|
L_{BA} {\mathsf P}_1^D
 {\mathsf e}_D(k) \otimes \left \langle {\mathsf e}_D(k) \right|$ is regarded as 
the Fourier transformation of the 
``resonance of cross species mass diffusion''. 
\end{remark} 

Now we will use Lemma \ref{lemmaba} (ii). Let $\Phi$ deonte any linear combination of $1,\xi_i,|\xi|^2$. We have that the operator
$$L_{AB}{\mathsf J}+ \frac{\sqrt{{\mathsf M}_B}}{  \sqrt{{\mathsf M}_A}  } L_{BA} {\mathsf J} $$
is orthogonal to all collision invariants $\sqrt{{\mathsf M}_A}, \xi_i\sqrt{{\mathsf M}_A}, |\xi|^2\sqrt{{\mathsf M}_A}$.
With this applied to ${\mathsf J} = L_{AB}^{-1}\xi_1{\mathsf E}_D$ (where we recall that ${\mathsf E}_D = {\mathsf M}_A$) it follows 
$$\int \left(L_{AB}(L_{AB}^{-1}\xi_1 {\mathsf E}_D) + \frac{\sqrt{{\mathsf M}_B}}{  \sqrt{{\mathsf M}_A}  } L_{BA} (L_{AB}^{-1}\xi_1 {\mathsf E}_D) \right) \Phi \sqrt{\mathsf M}_A d\xi =0$$
so that
$$\int \Phi\sqrt{{\mathsf M}_B} L_{BA} (L_{AB}^{-1}\xi_1 {\mathsf E}_D)d\xi = -\int \Phi \xi_1 {\mathsf M}_A d\xi$$
or so in the case of our interest

$$ {\mathsf E}_2 \otimes \left \langle {\mathsf E}_2 \right| L_{BA} L^{-1}_{AB} {\mathsf P}_1^D \xi^1 {\mathsf E}_D = $$
\begin{equation}
 \label{53}
 =  {\mathsf E}_2 \int_{\mathbb{R}^3_\xi}
 ( -\sqrt{3/2}\chi_0+\chi_4  ) L_{BA} L^{-1}_{AB} \xi^1 {\mathsf E}_D d\xi= - {\mathsf E}_2 \int_{\mathbb{R}^3_\xi}
 \frac{( -\sqrt{3/2}\chi_0+\chi_4  )}{\sqrt{{\mathsf M}_B}  }\xi^1 {\mathsf M}_A d\xi = 0.
\end{equation}
since the integrand is odd. This is an expression of the conservation of mass, energy and momentum in the entire system.
\begin{remark}
We could have proved the above also using the reflection preserving properties of the kernels and the collision frequency which is implied by Lemmas \ref{lemma32} and \ref{lemmaba} for both the operators $L_{AB}$ and $L_{BA}$ Considering the transformation $\xi \rightarrow -\xi$ due to the properties 
$$\nu^{AB}(-\xi)=\nu^{AB}(\xi) \text{, and  } k^{AB}(-\xi,\xi_*)=k^{AB}(\xi,-\xi_*)$$
if $f(-\xi)=-f(\xi)$ we get 
$$L^{AB}f(-\xi)=\int k^{AB}(-\xi,\xi_*) f(\xi_*)d\xi_* - \nu(-\xi)f(-\xi)=\int k^{AB}(\xi,-\xi_*)
 f(\xi_*)d\xi_* + \nu(\xi)f(\xi) = -L^{AB}f(\xi),$$ 
and similarly for $L_{BA}$. And similarly for an even function. From this it follows that the operator $L_{BA}L_{AB}^{-1}$ preserves the oddness of a function in $\xi$.
\end{remark}
From  \eqref{52} and \eqref{53}
\begin{multline} \label{54}
\sum_{j=1}^3 
e^{\sigma_j(k)(t-s) +\lambda(k) s}
 {\mathsf e}_j(k) \otimes \left \langle {\mathsf e}_j(k) \right|
L_{BA} {\mathsf P}_1^D
 {\mathsf e}_D(k) \otimes \left \langle {\mathsf e}_D(k) \right|
\\ = 
\sum_{j \in \{1,3\}}
e^{\sigma_j(k)(t-s) +\lambda(k) s}  \left\{
 -ik  {\mathsf E}_j \otimes \left \langle {\mathsf E}_j \right|
L_{BA} {\mathsf P}_1^D
 {\mathsf E}_D \otimes \left \langle {\mathsf E}_D \right| +O(1) k^2 \right\} 
+ e^{\sigma_2(k)(t-s) + \lambda(k)s} O(1) k^2.
\end{multline} 

Now from \eqref{54},  \eqref{spec_d2}, and  \eqref{spec_d3} we can write after performing the time integration 
\begin{multline} \label{as2}
\hat{\mathsf q}(x,t)=\int_{t^{\frac{1}{2} }}^{t-{t}^{\frac{1}{2} }} e^{(-i k \xi^1 +L )(t-s)}   L_{BA} e^{(-i k \xi^1 +L_{AB} )s} ds =\\  =
\sum_{j \in \{1,3\}}
\frac{ -ik \bigl( e^{\lambda(k) (t -t^{\frac{1}{2}}) + \sigma_j(k) t^{\frac{1}{2}} } - 
e^{\sigma_j(k)(t- t^{\frac{1}{2}}) + \lambda(k) t^{\frac{1}{2} }}
\bigr)  }{  \lambda(k) -\sigma_j(k)  } {\mathscr O}_j(k)
+  e^{\sigma_2(k)t } \frac{ k^2 
\bigl( e^{(\lambda(k) - \sigma_2(k))(t-t^{\frac{1}{2}}) } -1
\bigr)}{ (\sigma_2(k) - \lambda(k))    } {\mathscr O}_2(k) 
+ {\mathscr H}(k,t),
\end{multline}
where ${\mathscr O}_j$, $j=1,2,3$, are analytic functions around $k=0$ and independent of $t$.   

The function 
${\mathscr H}(k,t)$ is an analytic function given by
\begin{multline}
{\mathscr H}(k,t) \equiv
\int_{t^{\frac{1}{2}}}^{t- t^{\frac{1}{2}}} \sum_{j=1}^3 
e^{\sigma_j(k) (t-s) } {\mathsf e}_j(k) \otimes \left \langle {\mathsf e}_j(k) \right| 
 L_{BA} {\mathsf P}_1^D e^{ ( -i k \xi^1 +L_{AB})s} 
 \Pi_k^{D\perp} {\mathsf g}_{in} ds \\
+ \int_{t^{\frac{1}{2}}}^{t- t^{\frac{1}{2}}}  
e^{ ( -i k \xi^1 +L)(t-s)} 
\Pi_k^{\perp}
 L_{BA} {\mathsf P}_1^D 
e^{\lambda (k) s } {\mathsf e}_D(k) \otimes \left \langle {\mathsf e}_D(k) \right| 
 {\mathsf g}_{in} ds \\
+ \int_{t^{\frac{1}{2}}}^{t- t^{\frac{1}{2}}} 
 e^{ ( -i k \xi^1 +L)(t-s)} 
\Pi_k^{\perp}
 L_{BA} {\mathsf P}_1^D e^{ ( -i k \xi^1 +L)s} 
\Pi_k^{D\perp} {\mathsf g}_{in} ds.
\end{multline}
Due to the spectral properties \eqref{spec_d2} and \eqref{spec_d3}, one has that 
for $|k| \le \kappa_0$
\begin{equation}
\nonumber 
\| {\mathscr H}(k,t)\|_{L^2_\xi} \le O(1) e^{-t^{1/2}/C} \text{ for some } C>0.
\end{equation} 
This yields that 
\begin{equation}
\label{as3_1}
\left\| \int_{|k| \le \kappa_0} e^{ikx} {\mathscr H}(k,t)dk \right \|_{L^\infty_x(L^2_\xi)} \le O(1) e^{-t^{\frac{1}{2}}/C}.
\end{equation}

From the analytic property of the eigenvalues $\lambda(k)$ and  $\sigma_j(k)$ and the expansion properties in \eqref{spin1} and \eqref{B1.5}, one has that for $j=1,3$
\begin{equation} \label{as4} 
\displaystyle  \frac{ -ik }{  \lambda(k) -\sigma_j(k)  } {\mathscr O}_j(k) = 
\frac{ -ik }{  -a_2k^2  + i\lambda_jk + O(1)k^2} {\mathscr O}_j(k) \\
\end{equation}
is an analytic function of $k$ for $k<\kappa_0$ which is independent of t. 
>From \eqref{as4}, by shifting the contour of integration using the analyticity as done in \cite{liuyug1}, one has that for $|x| \le 2 |\lambda_1|t$
\begin{multline}
\label{as5}
\begin{cases}
\displaystyle  \left\|
\int_{|k| \le \kappa_0} e^{ikx} \frac{ k  e^{\lambda(k) (t -t^{\frac{1}{2}})
 + \sigma_j(k) t^{\frac{1}{2}}} }{  \lambda(k) -\sigma_j(k)  } {\mathscr O}_j(k) dk \right \|_{L^2_\xi} 
\le O(1) \frac{ e^{ -\frac{ x^2  }{  C(t+1)  }  } }{  \sqrt{(1+t)}  }  
\text{ for } j=1,3
\\
\displaystyle  \left\|
\int_{|k| \le \kappa_0} e^{ikx} \frac{ k 
e^{\sigma_j(k)(t- t^{\frac{1}{2}}) + \lambda(k) t^{\frac{1}{2} }}
  }{  \lambda(k) -\sigma_j(k)  } {\mathscr O}_j(k) dk \right \|_{L^2_\xi} 
\le O(1) \frac{ e^{ -\frac{ (x- \lambda_j t)^2  }{  C(t+1)  }  } }{  \sqrt{(1+t)}  } \text{ for } j=1,3
\end{cases}
\end{multline}
We just mention that the new contour of integration is of the form $\Omega$=$\Gamma_1 \cup \Gamma_2 \cup \Gamma_3$ where $\Gamma_1=(-\kappa_0/2, -\kappa_0/2+ic)$, $\Gamma_2=(-\kappa_0/2+ic,\kappa_0/2+ic)$, $\Gamma_3=(\kappa_0/2+ic,\kappa_0/2)$. Here the constant $c<\kappa_0$ is appropriately chosen as in \cite{liuyug1}.

Next we can write
\begin{multline} 
\displaystyle
e^{\sigma_2(k)t } \frac{ k^2  \bigl( e^{(\lambda(k) - \sigma_2(k))(t-t^{\frac{1}{2}}) } -1 \bigr) }{ (\sigma_2(k) - \lambda(k))    } {\mathscr O}_2(k) = \\ 
= e^{( -\frac{ 1 }{2} a_2k^2+O(k^3))t}  k^2(t-t^{\frac{1 }{2} }) e^{(-\frac{ 1 }{2} a_2k^2+O(k^3))t }\int _0^1 e^{(\lambda(k) - \sigma_2(k))(t-t^{\frac{1}{2}})u} du {\mathscr O}_2(k)=\\
=e^{ (-\frac{ 1 }{2} a_2k^2+O(k^3))t}  k^2(t-t^{\frac{1 }{2} }) {\mathscr O}^a_2(k,t)
\end{multline}
where ${\mathscr O}^a_2(k,t)$ is clearly an analytic function of $k$ and can be easily seen to be bounded in $t$ along the contour $\Omega$ by a constant independent of $k$ and $x$ provided $k$ is sufficiently small and $|x|<|2\lambda_1t|$. This allows one to again shift the contour of integration as in \cite{liuyug1} and obtain that for $|x| \le 2 |\lambda_1|t$

\begin{equation} \label{as51}
\displaystyle \left\| \int_{|k| \le \kappa_0} e^{ikx}   e^{\sigma_2(k)t } \frac{ k^2 
\bigl( e^{(\lambda(k) - \sigma_2(k))(t-t^{\frac{1}{2}}) } -1
\bigr)}{ (\sigma_2(k) - \lambda(k))    } {\mathscr O}_2(k) d k  \right \|_{L^2_\xi} 
\le t  O(1) \frac{ e^{ -\frac{ x^2  }{  C(t+1)  }  } }{  (1+t)^{3/2}  } = 
O(1) \frac{ e^{ -\frac{ x^2  }{  C(t+1)  }  } }{  (1+t)^{1/2}  }
\end{equation}

Finally, from the estimates in \eqref{as3_1},  \eqref{as5}, and \eqref{as51} combined with the spectral decomposition of $\hat {\mathsf q}(k)$ in \eqref{as2} and the definition \eqref{two-dec_d} one has that 
\begin{equation}
\label{as6}
\| {\mathsf q}_L(x,t) \|_{L^2_\xi}  \le O(1) \left( 
\sum_{j=1}^3 
\frac{ e^{ - \frac{ ( x- \lambda_j t)^2  }{  C(1+t)  }   } }{  \sqrt{(1+t)}  } 
+ e^{- t^{\frac{1}{2} }/C } 
 \right) \text{ for } |x| \le 2 |\lambda_1|(t+1).
\end{equation} 
Now from \eqref{as6}, \eqref{cons-tw2}, \eqref{51c}, and \eqref{51d} we conclude the estimate in the region $|x| \le 2 |\lambda_1 t|$:

\begin{equation} \label{ans_aa}
\|{\mathsf h}(x,t)\|_{L^2_\xi} \le O(1) \left( \sum_{j=1}^3 \frac{ e^{ - \frac{ ( x- \lambda_j t)^2  }{  C(1+t)  }   } }{  \sqrt{(1+t)}  } 
+ e^{- t^{\frac{1}{2} }/C } 
 \right)
\end{equation}

\noindent {\bf 2)  Estimate in the region $|x| \ge 2 |\lambda_1 t|$} 

 We simply use Duhamel's principle, 
\begin{multline}
\label{final_ans}
\|{\mathsf h}(x,t)\|_{L^2_\xi} 
\le \left\| \int_{{\mathbb R}} {\mathbb G}(x-y,t) {\mathsf h}_{in}(y) dy \right\|_{L^2_\xi} \\ + \left \|
 \int_0^t \int_{{\mathbb R}^2} 
{\mathbb G}(x-y,t-s)   
L_{BA} {\mathbb G}_{AB}(y-z,s) {\mathsf g}_{in}(z)
dz dy ds  \right\|_{L^2_\xi} \\
\le 
O(1) \left( 
\sum_{j=1}^3 
\frac{ e^{ - \frac{ ( x- \lambda_j t)^2  }{  C(1+t)  }   } }{  \sqrt{(1+t)}  } 
+ e^{- (|x|+t)/C } 
 \right)  \\ +
O(1) \int_0^t \int_{\mathbb R} \left(
\sum_{j=1}^3 
\frac{ e^{ - \frac{ ( x-y -\lambda_j (t-s))^2  }{  C(1+t)  }   } }{  \sqrt{(1+(t-s))}  } 
+ e^{- (|x-y|+(t-s))/C } 
 \right) \left(
\sum_{j=1}^3 
\frac{ e^{ - \frac{ ( y- \lambda_j s)^2  }{  C(1+t)  }   } }{  \sqrt{(1+s)}  } 
+ e^{- (|y|+s)/C } 
 \right)
dy ds \\ 
 \le 
O(1) \left( 
\sum_{j=1}^3 
\frac{ e^{ - \frac{ ( x- \lambda_j t)^2  }{ 2 C(1+t)  }   } }{  \sqrt{(1+t)}  } 
+ e^{- (|x|+t)/(2C) } 
 \right) \text{ for } |x| \ge 2 |\lambda_1 |(t+1).
\end{multline}

Finally, \eqref{basictheoa}, 
\eqref{ans_aa},  \eqref{final_ans} conclude Theorem \ref{maintheolin}.

\end{section}


\begin{thebibliography}{99}
\bibitem{Takata} Aoki, K.;  Bardos, C.; Takata, S. \textit{Knudsen layer for gas mixtures.}  J. Statist.
 Phys.  {\bf 112}  (2003),  no. 3-4, 629--655.


\bibitem{Kato} Kato, T. \emph{Perturbation theory for linear operators}, Die Grundlehren der mathematischen
 Wissenscheften, Band 132 Springer-Verlag New York, Inc., New York 1966.

\bibitem{Nik} Nicolaenko, B. \emph{Dispersion Laws for plane wave propagation}, In \emph{The Boltzmann 
Equation}, ed. by, F. A. Grunbaum, Courant Institute of Mathematical Sciences, 1971.



\bibitem{boltzmann} Boltzmann, L. (translated by Stephen G. Brush), 
\textit{Lectures on Gas Theory}, Dover Publications, Inc. New York, 
1964 

\bibitem{carleman} Carleman, T. Sur La Th\'eorie de l'\'Equation 
Int\'egrodiff\'erentielle de Boltzmann. \textit{ Acta Mathematica} {\bf 60} (1933), 91-142. 


\bibitem{Cerc} Cercignani, C.; Illner, R.; Pulvirenti, M. \emph{The Mathematical Theory of Dilute Gases},
 Applied Mathematical Sciences, 106. Springer-Verlag, New York, 1994. 
 


\bibitem{pinsky} Ellis, R.; Pinsky, M. The first and second 
fluid approximations to the linearized Boltzmann  equation. \textit{ J. Math. 
Pures Appl.} {\bf 54} (1975), no. 9, 125--156.  

\bibitem{bird} Hirschfelder, J.O.; Curtiss, C. F.; Bird, R. B., \textit{Molecular theory of Gases and Liquids}, John Wiley and Sons, New York, 1954
 
 
 
\bibitem{hilbert} Hilbert, D. \textit{Grundz\"uge einer Allgemeinen Theorie 
der Linearen Integralgleichungen}, (Teubner, Leipzig), Chap. 22 

 
\bibitem{kato} Kato, T. \textit{Perturbation theory for linear operators.}
 Die Grundlehren der mathematischen Wissenschaften, Band 132 Springer-Verlag
 New York, Inc., New York, 1966  
\bibitem{landau}     Lifshitz, E.M.;   Pitaevskii, L.P. \textit{Physical Kinetics},  Volume 10 (Course of Theoretical Physics), Butterworth-Heinemann; Reprint edition  (1981)

 

\bibitem{liuyu} Liu, T.-P.; Yu, S.-H. Boltzmann  Equation: 
Micro-Macro Decompositions  and Positivity of Shock Profiles. \textit{ Comm. 
Math. Phys.}, {\bf 246} (2004), no. 1, 133-179

\bibitem{liuyug1} Liu, T.-P.; Yu, S.-H.
The Green's Function and large-Time Behavior of
  Solutions for One-Dimensional Boltzmann Equation.
\textit{Comm. Pure Appl. Math.} {\bf 57} (2004), 1543-1608. 

\bibitem{liu-zeng} Liu, T.-P.; Zeng, Y. \textit{Large time behavior of 
solutions for general quasilinear  hyperbolic-parabolic systems of 
conservation laws}. Mem. Amer. Math. Soc. {\bf 125} (1997), no. 599,  
 


\bibitem{sone} Sone, Y. Kinetic Theory and Fluid Dynamics. \textit{Birkhauser} 2002. 


\bibitem{takata3} S. Takata; \textit{Kinetic theory analysis of the two-surface problem of a vapor-vapor
mixture in the continuum limit}, Phys. Fluids, {\bf 16}, no. 16, 2182-2188 (2004)



\bibitem{takata2} S. Takata; K. Aoki, \textit{The Ghost Effect in 
the continuum limit for a vapor gas mixture around condensed phases: 
Asymptotic analysis of the Boltzmann equation},
Transport Theory Statist. Phys. {\bf 30}, 205-237 (2001) 




 


 
\end{thebibliography}
\end{document}